\documentclass[conf]{new-aiaa}
\usepackage[utf8]{inputenc}
\usepackage{textcomp}

\usepackage{graphicx}
\usepackage{amsmath}
\usepackage[version=4]{mhchem}
\usepackage{siunitx}
\usepackage{longtable,tabularx}
\usepackage{algorithm, algorithmicx, algpseudocode}
\usepackage[acronym,nonumberlist]{glossaries-extra}
\usepackage[capitalize]{cleveref}
\usepackage{booktabs}
\usepackage{subcaption}
\usepackage{graphicx}
\usepackage{float}
\usepackage{mathtools}
\usepackage{multicol}
\usepackage{multirow}
\usepackage{physics}
\usepackage{silence}
\usepackage{siunitx}
\usepackage{ulem} 

\hypersetup{colorlinks=true,allcolors=cyan}

\newacronym{ads}{ADS}{Automatic Domain Splitting}
\newacronym{aop}{AOP}{Argument of Periapsis}

\newacronym{cam}{CAM}{Collision Avoidance Maneuver}
\newacronym{ci}{CI}{Confidence Interval}
\newacronym{cut}{CUT}{Conjugate Unscented Transform}
\newacronym{cut4}{CUT4}{Fourth-Order Conjugate Unscented Transform}
\newacronym{cut6}{CUT6}{Sixth-Order Conjugate Unscented Transform}

\newacronym{da}{DA}{Differential Algebra}
\newacronym{dmm}{DMM}{Dirac Mixture Model}
\newacronym{dod}{DoD}{Department of Defense}

\newacronym{ekf}{EKF}{Extended Kalman Filter}
\newacronym[\glslongpluralkey={Equations of Motion}]{eom}{EOM}{Equation of Motion}

\newacronym{fp}{FP}{Floating Point}

\newacronym{gcrf}{GCRF}{Geocentric Celestial Reference Frame}
\newacronym{gmm}{GMM}{Gaussian Mixture Model}

\newacronym{heo}{HEO}{High Earth Orbit} 
\newacronym{hf}{HF}{high-fidelity}
\newacronym{hsl}{HSL}{Harwell Subroutine Library}

\newacronym{ic}{IC}{Initial Condition}
\newacronym{ipopt}{IPOPT}{Interior Point Optimizer}

\newacronym{jpl}{JPL}{Jet Propulsion Laboratory}

\newacronym{lam}{LAM}{Likelihood Agreement Measure}
\newacronym{lc}{LinCov}{Linear Covariance}
\newacronym{leo}{LEO}{Low Earth Orbit}
\newacronym{lf}{LF}{low-fidelity}
\newacronym{loads}{LOADS}{Low-Order Automatic Domain Splitting}

\newacronym{mce}{MCE}{Moment Constraint Equation}
\newacronym{mc}{MC}{Monte Carlo}
\newacronym{mee}{MEE}{Modified Equinoctial Element}
\newacronym{meo}{MEO}{Medium Earth Orbit}
\newacronym{mf}{MF}{multifidelity}

\newacronym{nli}{NLI}{Nonlinearity Index}
\newacronym{nlp}{NLP}{Nonlinear Programming Problem}
\newacronym{norad}{NORAD}{North American Aerospace Defense Command}

\newacronym{od}{OD}{Orbit Determination}
\newacronym{ode}{ODE}{Ordinary Differential Equation}
\newacronym{oe}{OE}{Orbital Element}
\newacronym[longplural={Orders of Magnitude}]{oom}{OOM}{Order of Magnitude}

\newacronym{pce}{PCE}{Polynomial Chaos Expansion}
\newacronym{pdf}{PDF}{Probability Density Function}
\newacronym{pw}{PW}{point-wise}

\newacronym{raan}{RAAN}{Right Ascension of the Ascending Node}
\newacronym{rhs}{RHS}{Right-Hand Side}
\newacronym{rmse}{RMSE}{Root Mean Square Error}

\newacronym{sdp}{SDP}{Simplified Deep Space Perturbations}
\newacronym{sgp}{SGP}{Simplified General Perturbations}
\newacronym{srp}{SRP}{Solar Radiation Pressure}
\newacronym{ssa}{SSA}{Space Situational Awareness}
\newacronym{sst}{SST}{Space Surveillance and Tracking}
\newacronym{stm}{STM}{State Transition Matrix}
\newacronym{stt}{STT}{State Transition Tensor}

\newacronym{tca}{TCA}{Time of Closest Approach}
\newacronym{teme}{TEME}{True Equator Mean Equinox}
\newacronym{tle}{TLE}{Two-Lines Element Set}

\newacronym{ukf}{UKF}{Unscented Kalman Filter}
\newacronym{up}{UP}{Uncertainty Propagation}
\newacronym{us}{US}{United States}
\newacronym{ut}{UT}{Unscented Transform}
\newacronym{utc}{UTC}{Coordinated Universal Time}

\glsdisablehyper

\newcolumntype{C}{>{$}c<{$}} 
\newcolumntype{R}{>{$}r<{$}} 
\newcolumntype{L}{>{$}l<{$}} 

\algnotext{EndFunction}

\DeclareMathOperator*{\argmax}{arg\,max}

\title{Multifidelity Orbit Uncertainty Propagation using\\Taylor Polynomials}

\author{Alberto Foss\`a\footnote{PhD Candidate, Department of Aerospace Vehicles Design and Control, 10 Avenue Edouard Belin.}}
\affil{Institut Sup\'erieur de l’A\'eronautique et de l’Espace, Toulouse, 31055, France}

\author{Roberto Armellin\footnote{Professor of Astrodynamics, Te P\=unaha \=Atea - Space Institute, 20 Symonds Street.}}
\affil{The University of Auckland, Auckland, 1010, New Zealand}

\author{Emmanuel Delande\footnote{Space Surveillance Engineer, Space Surveillance and Tracking, 18 Avenue Edouard Belin.}}
\affil{Centre National d'\'Etudes Spatiales, Toulouse, 31401, France}

\author{Matteo Losacco\footnote{Postdoctoral Researcher, Department of Aerospace Vehicles Design and Control, 10 Avenue Edouard Belin.}}
\affil{Institut Sup\'erieur de l’A\'eronautique et de l’Espace, Toulouse, 31055, France}

\author{Francesco Sanfedino\footnote{Associate Professor, Department of Aerospace Vehicles Design and Control, 10 Avenue Edouard Belin.}}
\affil{Institut Sup\'erieur de l’A\'eronautique et de l’Espace, Toulouse, 31055, France}

\begin{document}

\maketitle

\begin{abstract}
A new multifidelity method is developed for nonlinear orbit uncertainty propagation. This approach guarantees improved computational efficiency and limited accuracy losses compared to fully high-fidelity counterparts. The initial uncertainty is modeled as a weighted sum of Gaussian distributions whose number is adapted online to satisfy the required accuracy. As needed, univariate splitting libraries are used to split the mixture components along the direction of maximum nonlinearity. Differential Algebraic techniques are used to propagate these Gaussian kernels and compute a measure of nonlinearity required for the split decision and direction identification. Taylor expansions of the flow of the dynamics are computed using a low-fidelity dynamical model to maximize computational efficiency and corrected with selected high-fidelity samples to minimize accuracy losses. The effectiveness of the proposed method is demonstrated for different dynamical regimes combining SGP4 theory and numerical propagation as low- and high-fidelity models respectively.
\end{abstract}

\section{Introduction}\label{sec:introduction}

\lettrine{A}{ccurate} propagation of uncertainties is key for \acrfull{ssa} applications. Ranging from \acrfull{sst} to \acrfull{cam} planning, these activities  require the estimation of the spacecraft state at a future epoch given the \acrfullpl{ic} at an earlier time. These \acrshortpl{ic} are the result of ground-based observations or on-board \acrfull{od} techniques, which are stochastic in nature. The propagation of a nominal state is thus not sufficient, but it is rather required to propagate the associated initial \acrfull{pdf} resulting from the estimation procedure. Given the intrinsic nonlinearity of the dynamics, no analytical solution exists to map the aforementioned \acrshort{pdf}. On the contrary, several methods exist to approximate this quantity with increasing level of accuracy, usually at the expenses of higher computational efforts. Instead, this paper focuses on a novel \acrlong{mf} approach to tackle the \acrfull{up} problem in orbital dynamics which guarantees a similar accuracy of \acrlong{hf} techniques with a reduced computational load.\par

Previous approaches tackled the \acrshort{up} problem in different ways, the most widely used being linear methods and sample-based techniques \cite{Luo2017}. Linear methods build on the assumptions that a linearized dynamical model is sufficient to capture the dynamics of neighboring trajectories and that the uncertainty can be completely characterized by a multivariate Gaussian distribution. Under these assumptions, only the mean and the covariance matrix need to be propagated, and the problem reduces to the integration of the \acrfull{stm} known as \acrfull{lc} propagation. However, if the dynamics is highly nonlinear, this approach fails to accurately describe the time evolution of the state \acrshort{pdf}. An extension of the \acrshort{lc} method is that of \acrfullpl{stt}, which are higher-order Taylor expansions of the dynamics about the nominal trajectory \cite{Park2006a}. However, this method is not ideal for complex dynamics at high orders since it requires the derivation of analytical expressions for the  variational equations and their numerical integration along with the reference trajectory.\par

On the other side, \acrfull{mc} simulations are widely used in operational scenarios since they provide highly accurate estimation of the state uncertainty and are easy to implement. Yet their accuracy comes with a major computational cost, which makes \acrshort{mc} techniques not applicable for real-time applications or maintenance of very large space objects catalogs. Sample based methods also include \acrfull{ut} and \acrfull{cut} \cite{Julier2004,Adurthi2015,Adurthi2018}, which are deterministic in nature as opposed to stochastic \acrshort{mc} simulations. These methods are based on the idea that is easier to approximate a state \acrshort{pdf} than an arbitrary nonlinear transformation. Assuming a Gaussian distribution for the initial state uncertainty, they enforce the corresponding \acrfullpl{mce} up to a given order to solve for a deterministic set of samples, denoted as sigma points, whose weighted statistics correctly capture the first statistical moments of the state \acrshort{pdf}. To reduce the search space, these points are constrained to lie on carefully selected set of axes, and a fully symmetric set is chosen to automatically satisfy the \acrshortpl{mce} for any odd moment of the Gaussian distribution. If the \acrshort{ut} can only capture the statistical moments up to third order, \acrshort{cut} sigma points sets were obtained to match these moments up to ninth order. However, this approach suffers from a rapid increase in the number of samples needed to satisfy the \acrshortpl{mce} for increasing orders and state dimension.\par

\acrfull{pce} is another samples-based technique developed for \acrshort{up} \cite{Jones2018}. This method aims at obtaining a functional representation of the propagated uncertainty with respect to the input random variables, thus providing accurate information on higher-order statistical moments of the state \acrshort{pdf}. Inputs and outputs of the considered transformation are represented by series approximations of standard random variables rewritten as series expansions of orthogonal polynomials for which the coefficients are sought. This technique requires a number of \acrshort{pce} terms which varies exponentially with both the order of the polynomial basis and the dimensions of the input random variables, leading to the course of dimensionality for large problems.\par

A new \acrlong{mf} technique based on Taylor polynomials is developed in this work to address the need for accurate and less computationally expensive methods for nonlinear \acrlong{up} in orbital dynamics. In this context, \acrlong{mf} techniques combine different fidelity force models to accurately propagate the state \acrshort{pdf} at a reduced cost \cite{Jones2018,Peherstorfer2016}. With sample-based methods such as \acrshort{pce}, this approach requires the propagation of most samples in \acrlong{lf} dynamics and the evaluation of the \acrlong{hf} model on few carefully chosen points to correct the initial approximate solution. On the contrary, the proposed method relies on a Taylor expansion of the \acrlong{lf} dynamics for which only the constant part needs to be propagated with the \acrlong{hf} model.\par

Taylor algebra, or \acrfull{da}, provides the tools to compute the arbitrary-order derivatives of a function within a computer environment \cite{Berz1999}. In the context of \acrshort{up}, \acrshort{da} is similar to \acrshortpl{stt} in the sense that it provides an expansion of the flow of the dynamics around the nominal state in the form of a Taylor polynomial of order $n$. This information was successfully employed in previous work to map the initial state uncertainties through the nonlinear orbital dynamics \cite{Valli2013,armellin2018probabilistic}. However, if the transformation under consideration is strongly nonlinear, a single Taylor polynomial might fail to accurately describe the \acrshort{ode} flow in the domain of interest around the nominal solution. \acrfull{ads} techniques were developed to overcame this issue \cite{Wittig2015a}. The idea is to continuously monitor the accuracy of the Taylor expansion by estimating the magnitude of the truncated terms (terms of order $n + 1$) and react consequently when the imposed error threshold is violated. When it happens, the single polynomial is split into two new expansions, each of them covering only half of the initial domain, so as to meet the specified accuracy. This algorithm requires an expansion order $\geq 3$ and its accuracy increases for increasing orders. In this work, a novel splitting algorithm specifically developed for second-order Taylor expansions and named \glsentryfull{loads} is employed instead \cite{Losacco}.\par

\acrfullpl{gmm} can approximate an arbitrary \acrshort{pdf} using a weighted sum of Gaussian distributions, and this approximation converges to the true \acrshort{pdf} as the number of components approaches infinity. In orbit \acrshort{up}, algorithms were developed to adapt online the number of \acrshort{gmm} components required to correctly capture the entire state \acrshort{pdf} as the assumption of Gaussianity becomes less accurate for longer propagation times or stronger nonlinearities \cite{DeMars2013}. In this scenario, the initial uncertainty is described using a single Gaussian kernel and a metric is used throughout the propagation step to continuously monitor the error between linearized and fully nonlinear models. When the specified threshold is violated, the propagation is halted to split the single kernel into components and then resumed on the newly obtained kernels. Other works \cite{Horwood2011,Vittaldev2016} simplify this approach by performing all required splits on the initial distribution and keep the number of Gaussian kernels fixed for the whole propagation time span. However, compared to their adaptive counterparts, these methods suffer from an increased computational effort due to the large number of kernels to be handled in the early stages of propagation. Hybrid methods were proposed that combine \acrshortpl{gmm} with \acrfull{ekf}, \acrfull{ukf} \cite{DeMars2013}, \acrshortpl{pce} \cite{Vittaldev2016}, and \acrshort{da} \cite{sun2019nonlinear}, but the Taylor expansions of the dynamics were never used to inform the splitting into the \acrshort{gmm} components.\par

In this context, a novel \acrfull{mf} method that combines \acrshort{da}-based propagation and adaptive splitting of \acrshortpl{gmm} in \acrfull{lf} dynamics with \acrfull{pw} propagation of the polynomials centers in \acrfull{hf} is presented. In the \acrshort{lf} step, a second-order Taylor expansion of the flow is used to compute the \acrlong{nli} $\nu$. When its magnitude crosses a given threshold, the polynomial is split into three new expansions such that the assumption of a quasi-linear transformation still holds in each subdomain. Assuming a Gaussian distribution for the initial state, the \acrshort{ic} is constructed such that the former distribution is correctly captured within the domain of validity of its Taylor polynomial. As the propagation steps forward, more domains are generated by the \acrshort{loads} algorithm, each of them corresponding to a new \acrshort{gmm} component such that their weighted sum correctly captures the true non-Gaussian state \acrshort{pdf}. Once the \acrshort{da}-based propagation has ended, evaluation of the resulting polynomials on \acrshort{ut} sigma-points is used to map the statistics of the initial components to final time. The \acrshort{hf} step is then simply a \acrshort{pw} propagation of the means of these Gaussian kernels for the same time span. The resulting final states will constitute the polynomials centers and components means of the updated \acrshort{mf} solution.\par

Nonlinearity of any dynamical system is not an intrinsic attribute of the system itself, but is rather coupled with the mathematical representation chosen to describe its geometry, kinematics and evolution dynamics \cite{Junkins2004}. In the context of orbit \acrshort{up}, this aspect is emphasized by the choice of \acrfullpl{oe} used to describe the osculating orbit state. As a consequence, the performance of the developed \acrshort{loads}-\acrshort{gmm} algorithm is expected to be strongly affected by this choice of state parametrization, being the splitting decision derived from a measure of nonlinearity of the neighboring dynamics. If cartesian coordinates are widely used in end-user applications, their evolution over time is strongly nonlinear due to the Keplerian attraction being proportional to $1/r^2$. Several \acrshortpl{oe} sets adapted to the underlying dynamics were thus proposed to mitigate these nonlinearities. Examples include Keplerian elements, equinoctial elements \cite{Broucke1972} and \acrlongpl{mee} \cite{Walker1985}. Under the assumption of unperturbed Keplerian motion, five of these coordinates describe the shape and spatial orientation of the orbit and are thus constant in time, while the last provides the instantaneous position of the object along its closed path. More recently, \citet{Horwood2011} introduced a set of alternate equinoctial elements which, in absence of perturbations, evolve linearly in time demonstrating their effectiveness for orbit \acrshort{up} using Gaussian sum filters. The hypothesis of restricted two-body dynamics was then removed by \citet{Aristoff2021a} who proposed a set of $J_2$ equinoctial orbital elements that absorb the perturbation due to the equatorial bulge of the central attracting body. Finally, \citet{Bau2021} introduced the generalized equinoctial orbital elements which extend this approach to any perturbing force derived from a potential.

The paper is organized as follows. After a short description of the \acrshort{da} formalism in \cref{sec:uncertainty_propagation}, the \acrshort{loads} algorithm is described in \cref{sec:loads}. \acrshortpl{gmm} are then introduced in \cref{sec:gmm}. \Cref{sec:univariate_splitting_lib} describes the univariate splitting libraries for standard normal distributions while \cref{sec:multivariate_splitting} discusses their extension to the multivariate case. \Cref{sec:loads_gmm_algorithm} presents the combined \acrshort{loads}-\acrshort{gmm} algorithm for \acrshort{up} that constitutes the \acrshort{lf} step of the proposed method. The \acrshort{mf} approach is then described in \cref{sec:mf_orbit_up}. Three numerical applications to the problem of orbit \acrshort{up} are then included in \cref{sec:numerical_results} and the accuracy of the \acrshort{mf} method compared against reference \acrshort{mc} simulations. The importance of the \acrshort{hf} correction and the computational efficiency of the overall method are then assessed in \cref{sec:rmse,sec:computational_load} respectively. Conclusions are finally drawn in \cref{sec:conclusions}.\par

\section{\glsentrytitlecase{up}{long}}\label{sec:uncertainty_propagation}

Consider the generic nonlinear transformation $\vb*{f}\colon\mathbb{R}^n\to\mathbb{R}^m$

\begin{equation}
    \vb*{y} = \vb*{f}(\vb*{x})
    \label{eqn:generic_nonlinear_transform}
\end{equation}
and two multivariate random variables $\vb*{X},\vb*{Y}$ defined on $\mathbb{R}^n,\mathbb{R}^m$ with \acrshortpl{pdf} $p_{\vb*{X}},p_{\vb*{Y}}$ respectively. The problem of nonlinear \acrshort{up} consists in estimating $p_{\vb*{Y}}$ from the only knowledge of $p_{\vb*{X}}$ and $\vb*{f}$. In \acrshort{da} framework, $\vb*{X}$ is expressed as a function of its expected value and associated uncertainty as

\begin{equation}
    [\vb*{x}] = \bar{\vb*{x}} + \vb*{\beta} \cdot \delta\vb*{x}
\end{equation}
with $\bar{\vb*{x}}=\mathbb{E}[\vb*{X}]$, $\vb*{\beta}$ \acrfull{ci} semi-amplitude and $\delta\vb*{x}$ first-order variations around $\bar{\vb*{x}}$.

\clearpage
\noindent
Evaluating $\vb*{f}$ in the \acrshort{da} framework results in the Taylor expansion of $\vb*{y}$ in terms of $\delta\vb*{x}$ expressed as

\begin{equation}
    [\vb*{y}] = \vb*{f}([\vb*{y}]) = \mathcal{T}_{\vb*{y}}(\delta\vb*{x})
    \label{eqn:taylor_expansion_generic_y}
\end{equation} 
The polynomials in \cref{eqn:taylor_expansion_generic_y} can then be used to efficiently compute the propagated statistics as described in \cref{sec:gmm}.\par

\subsection{\glsentrytitlecase{loads}{long}}\label{sec:loads}

The \acrshort{loads} algorithm is a \acrshort{da}-based method that aims at solving the uncertainty propagation problem by splitting the uncertainty set into a manifold of second-order Taylor expansions \cite{Losacco}. The splitting decision and direction are determined by a proper measure of nonlinearity that is derived from the analysis of the Taylor expansion of the Jacobian of the transformation. This \acrfull{nli} $\nu$ is a reformulation in \acrshort{da} sense of the same parameter introduced by \citet{Junkins2004} and approximates the relative error between full nonlinear and linearized transformations at the edges of the domains on which the Taylor polynomials are defined.

Consider the generic transformation $\vb*{f}$ and a second order expansion for $[\vb*{x}]$. The Jacobian of $\vb*{f}$ is defined as the $m \times n$ matrix $[\vb*{J}]$ of the partials of $\vb*{f}$ with respect to $[\vb*{x}]$

\begin{equation}
    [J_{ij}] = \dfrac{\partial f_i}{\partial x_j} \qquad i=1,\ldots,m\quad j=1,\ldots,n
    \label{eqn:generic_jacobian}
\end{equation}
If the transformation is linear, the Taylor expansion of the Jacobian coincides with its constant part, i.e. $[\vb*{J}] = \overline{\vb*{J}}$. On the contrary, if the transformation is nonlinear, each component $[J_{ij}]$ of $[\vb*{J}]$ is a first-order Taylor expansion in the deviations $\delta\vb*{x}$, i.e.

\begin{equation}
    [J_{ij}] = \overline{J}_{ij}+\delta J_{ij}=\mathcal{T}_{J_{ij}}(\delta\vb*{x})
    \label{eqn:jacobian_component}
\end{equation}
Consider now the non-constant part of the Jacobian matrix elements, $\delta J_{ij}$. Each element can be expressed as a linear combination of the variables $\delta x_k$ for $k=1,\ldots,n$, i.e.

\begin{equation}
    \delta J_{ij}=\mathcal{T}_{\delta J_{ij}}=\displaystyle\sum\limits_{k=1}^n a_{ij,k}\delta x_k
    \label{eqn:nonconstant_jacobian_component}
\end{equation}
as a result, a rigorous estimation of the variation range of each term $\delta J_{ij}$ can be obtained as

\begin{equation}
    \begin{aligned}
        b_{\delta J_{ij}} = \displaystyle\sum\limits_{k=1}^n \left\lvert a_{ij,k}\right\rvert\\
        \delta J_{ij}\in[-b_{\delta J_{ij}},b_{\delta J_{ij}}]
    \end{aligned}
    \label{eqn:jacobian_component_bounds}
\end{equation}
The idea behind the \acrshort{loads} algorithm is to use these bounds as a measure of nonlinearity of the transformation. More specifically, the \acrshort{nli} is expressed as

\begin{equation}
    \nu = \cfrac{\displaystyle\sum\limits_{i=1}^m\sum\limits_{j=1}^n b_{\delta J_{ij}}}{\displaystyle\sum\limits_{i=1}^m\sum\limits_{j=1}^n \left\lvert\overline{J}_{ij}\right\rvert}=\cfrac{\displaystyle\sum\limits_{i=1}^m\sum\limits_{j=1}^n\sum\limits_{k=1}^n \left\lvert a_{ij,k}\right\rvert}{\displaystyle\sum\limits_{i=1}^m\sum\limits_{j=1}^n \left\lvert\overline{J}_{ij}\right\rvert}
    \label{eqn:nonlinear_index_def}
\end{equation}
Therefore, the index $\nu$ is obtained as the ratio between the entry-wise 1-norm of the matrix of the bounds and the entry-wise 1-norm of the constant part of the Jacobian. In case of linear transformation, all elements $b_{\delta J_{ij}}$ are identically equal to zero, so $\nu=0$. As the transformation departs from linearity, the numerator increases, thus a hint on the nonlinear nature of the investigated problem is obtained.

In case nonlinearities are detected, the algorithm proceeds by splitting the uncertainty set along a determined direction $k$. The identification of the splitting direction is done with the same procedure. The only element that changes is the processed Jacobian. More specifically, given an investigated splitting direction $e$ and the associated variable $\delta x_e$, defining as $\delta\vb*{x}_e$ the vector

\begin{equation}
    \delta\vb*{x}_e = \left\{0,\ldots,0,\delta x_e,0,\ldots,0\right\}^{\textrm{T}}\qquad e\in[1,\ldots,n]
    \label{eqn:directional_deviation}
\end{equation}
a directional Jacobian is built as the composition of $[\vb*{J}(\delta\vb*{x})]$ with $\delta \vb*{x}_e$, i.e.

\begin{equation}
    [\left.\vb*{J}\right\rvert_e]=\mathcal{T}_{\vb*{J}}(\delta\vb*{x})\circ \delta\vb*{x}_e = \mathcal{T}_{\left.\vb*{J}\right\rvert_e}(\delta \vb*{x}_e)
    \label{eqn:directional_jacobian}
\end{equation}
The directional \acrshort{nli} is then computed as

\begin{equation}
    \nu_e = \cfrac{\displaystyle\sum\limits_{i=1}^m \sum\limits_{j=1}^n b_{\delta \left.J_{ij}\right\rvert_e}}{\displaystyle\sum\limits_{i=1}^m\sum\limits_{j=1}^n \left\lvert\overline{J}_{ij}\right\rvert}
    \label{eqn:directional_nonlinear_index_def}
\end{equation}
The splitting direction $k$ is finally selected as the one corresponding to the maximum directional index, i.e.

\begin{equation}
    k = \argmax_e\left\{\nu_e\right\}
    \label{eqn:splitting_direction}
\end{equation}
Once the splitting direction has been determined, the current set $[\vb*{x}]$ is split according to the preselected splitting library $\mathcal{S}_l$ as described in \cref{sec:loads_gmm_algorithm}, thus generating a list of $N_{sub}$ subsets $[\vb*{x}^{(p)}]$ with $N_{sub}$ size of the library itself.

The splitting algorithm is then embedded in an iterative procedure as described in \cref{alg:loads_algorithm} to generate the so-called manifold $M_{\vb*{y}}$ of subdomains $[\vb*{y}^{(p)}]$ that satisfy the imposed constraints. These constraints include a maximum nonlinearity threshold $\varepsilon_\nu$ and a maximum number of splits $N_{max}$ so that $\nu^{(p)}<\varepsilon_\nu \lor N_{splits}^{(p)}=N_{max}\ \forall\ [\vb*{y}^{(p)}] \in M_{\vb*{y}}$. The choice of $\varepsilon_\nu\in\mathbb{R}^+$ is thus a trade-off between accuracy of the \acrshort{da} maps, increasing for $\varepsilon_\nu\to 0$, and size of the output manifold $M_{\vb*{y}}$, decreasing for $\varepsilon_\nu\to +\infty$. However, being $\nu$ a direct measure of the relative error between full nonlinear and linearized transformations, $\varepsilon_\nu$ must be chosen sufficiently small ($\varepsilon_\nu\ll 1$) to guarantee an acceptable accuracy of the output \acrshort{da} maps. A more exhaustive analysis of \cref{alg:loads_algorithm} can be found in \cite{Losacco} where the \acrshort{loads} algorithm is applied to the investigated case of orbit \acrshort{up}, and a detailed analysis of the role played by the different control parameters is offered.

\begin{algorithm}[H]
\begin{algorithmic}
\setstretch{1.2}
    \Function{loads}{$\vb*{f},[\vb*{x}],\varepsilon_\nu,N_{max},\mathcal{S}_l$}
        \State Initialize the working and output manifolds $W_{\vb*{x}}$, $M_{\vb*{x}}, M_{\vb*{y}}$
        \State Add $[\vb*{x}]$ to $W_{\vb*{x}}$
        \While{$W_{\vb*{x}}$ is not empty}
            \State Remove the first set $[\vb*{x}^{(p)}]$ from $W_{\vb*{x}}$
            \State Evaluate $[\vb*{y}^{(p)}] = \vb*{f}([\vb*{x}^{(p)}])$
            \State Compute the \acrshort{nli} $v^{(p)}$ \Comment{see \cref{eqn:generic_jacobian,eqn:jacobian_component,eqn:nonconstant_jacobian_component,eqn:jacobian_component_bounds,eqn:nonlinear_index_def}}
            \If{$\nu^{(p)}<\varepsilon_{\nu} \lor N_{splits}^{(p)} = N_{max}$}
                \State Add $[\vb*{x}^{(p)}]$ to $M_{\vb*{x}}$ and $[\vb*{y}^{(p)}]$ to $M_{\vb*{y}}$
            \Else
                \State Compute the splitting direction $k$ \Comment{see \cref{eqn:directional_jacobian,eqn:directional_nonlinear_index_def,eqn:splitting_direction}}
                \State Split $[\vb*{x}^{(p)}]$ into $[\vb*{x}^{(i)}]$ for $i=1,\ldots,N_{sub}$ \Comment{see \cref{sec:loads_gmm_algorithm}}
                \State Add all $[\vb*{x}^{(i)}]$ to $W_{\vb*{x}}$
            \EndIf
        \EndWhile
        \State \Return $M_{\vb*{x}},M_{\vb*{y}}$
    \EndFunction
\end{algorithmic}
\caption{\glsentryshort{loads} algorithm}
\label{alg:loads_algorithm}
\end{algorithm}

\clearpage
\subsection{\glsentrytitlecase{gmm}{long}}\label{sec:gmm}

Consider a multivariate random variable $\vb*{X}$ defined on $\mathbb{R}^n$ with \acrshort{pdf} $p_{\vb*{X}}$. The probability density function $p_{\vb*{X}}$ has no general closed-form solution and might be approximated using different techniques. However, if $\vb*{X}$ is Gaussian with mean $\vb*{\mu}_{\vb*{X}}\in\mathbb{R}^n$ and covariance $\vb*{P}_{\vb*{X}}\in\mathbb{R}^{n \times n}$, its \acrshort{pdf} has the well-known analytical expression

\begin{equation}
    p_{\vb*{X}}(\vb*{x}) = \mathcal{N}(\vb*{x};\vb*{\mu}_{\vb*{X}},\vb*{P}_{\vb*{X}}) = \dfrac{1}{\sqrt{|2\pi\vb*{P}_{\vb*{X}}|}} \cdot \exp{-\frac{1}{2}(\vb*{x}-\vb*{\mu}_{\vb*{X}})^T \vb*{P}_{\vb*{X}}^{-1} (\vb*{x}-\vb*{\mu}_{\vb*{X}})}
    \label{eqn:multivariate_gaussian_pdf}
\end{equation}
where $|\cdot|$ represents the matrix determinant.

\acrshortpl{gmm} leverage this property to obtain analytic approximations of arbitrary \acrshortpl{pdf} as weighted sums of Gaussian distributions. For a generic random variable $\vb*{X}$, its \acrshort{pdf} is thus expressed as

\begin{equation}
    p_{\vb*{X}}(\vb*{x}) \simeq \sum_{i=1}^{L} \alpha_i p_{\mathcal{N}}(\vb*{x};\vb*{\mu}_i,\vb*{P}_i)
    \label{eqn:gmm_pdf_def}
\end{equation}
with $L$ total number of components and $\alpha_i$ mixing proportions or weights subject to

\begin{subequations}
    \begin{gather}
        \alpha_i \geq 0 \qquad \forall i \in \{1,\ldots,L\}\\
        \sum_{i=1}^{L} \alpha_i = 1
    \end{gather}
\end{subequations}
If $\vb*{X}$ is subject to a nonlinear transformation such as \cref{eqn:generic_nonlinear_transform}, the number of components $L$ that guarantees an accurate representation of the input \acrshort{pdf} $p_{\vb*{X}}$ might fail to correctly capture the shape of the transformed \acrshort{pdf} $p_{\vb*{Y}}$. A splitting algorithm is thus employed to dynamically adapt this number so that $p_{\vb*{Y}}$ is approximated with the desired accuracy. Its formulation is described in \cref{sec:univariate_splitting_lib,sec:multivariate_splitting,sec:loads_gmm_algorithm}.

\subsubsection{Univariate Splitting Library}\label{sec:univariate_splitting_lib}

Consider the univariate standard Gaussian distribution $X\sim\mathcal{N}(0,1)$ with \acrshort{pdf}

\begin{equation}
    p_X(x) = p_{\mathcal{N}}(x;0,1)
    \label{eqn:univariate_standard_normal_pdf}
\end{equation}
for the purpose of computing a univariate splitting library, it is desired to approximate \cref{eqn:univariate_standard_normal_pdf} with a \acrshort{gmm} such that

\begin{equation}
    p_X(x) \simeq \sum_{i=0}^{L-1} \tilde{\alpha}_i p_{\mathcal{N}}(x;\tilde{\mu}_i,\tilde{\sigma}_i^2)
    \label{eqn:univariate_gmm_def}
\end{equation}
with $L$ number of components and $\tilde{\alpha}_i,\tilde{\mu}_i,\tilde{\sigma}^2_i$ mixture weights, means and variances respectively. After selecting the number of components $L$, a minimization problem is formulated to determine the remaining $3L$ parameters. The number of free variables can be then further reduced constraining the \acrshort{gmm} to be homoscedastic, i.e.

\begin{equation}
    \tilde{\sigma}^2_i = \tilde{\sigma}^2 \ \forall i \in\{0,\ldots,L-1\}
\end{equation}
and taking advantage of the symmetric shape of $p_X$ about the mean. \Cref{eqn:univariate_gmm_def} is thus rewritten as

\begin{equation}
    p_X(x) \simeq \tilde{\alpha}_0 \cdot p_{\mathcal{N}}(x;0,\tilde{\sigma}^2) + \sum_{i=1}^{l} \tilde{\alpha}_i \cdot \bigg[p_{\mathcal{N}}(x;-\tilde{\mu}_i,\tilde{\sigma}^2) + p_{\mathcal{N}}(x;\tilde{\mu}_i,\tilde{\sigma}^2)\bigg]
    \label{eqn:univariate_gmm_symmetric}
\end{equation}
with $l=\lfloor L/2 \rfloor$ and $\alpha_0=0$ if $L$ is even so that only $L+1$ parameters remain to be determined.\par

To solve for the optimal $\tilde{\alpha_i},\tilde{\mu_i}$ and $\tilde{\sigma}^2$ a performance index $J$ is developed starting from the $L_2$ distance between \acrshortpl{pdf} defined as

\begin{equation}
    L_2(p_1,p_2) = \int_{\mathbb{R}^n} \left[p_1(\vb*{x}) - p_2(\vb*{x}) \right]^2 \mathrm{d}\vb*{x}
    \label{eqn:l2_pdf_def}
\end{equation}
The $L_2$ distance is preferred here instead of other measures since an analytical expression for \cref{eqn:l2_pdf_def} exists if both $p_1,p_2$ are approximated by \acrshortpl{gmm} as

\begin{subequations}
    \begin{gather}
        p_1(\vb*{x}) \simeq \sum_{i=1}^{L_1} \alpha_{1,i} p_{\mathcal{N}}(\vb*{x};\vb*{\mu}_{1,i},\vb*{P}_{1,i})\\
        p_2(\vb*{x}) \simeq \sum_{j=1}^{L_2} \alpha_{2,j} p_{\mathcal{N}}(\vb*{x};\vb*{\mu}_{2,j},\vb*{P}_{2,j})
    \end{gather}
    \label{eqn:p1_p2_gmm_def}
\end{subequations}
substituting \cref{eqn:p1_p2_gmm_def} into \cref{eqn:l2_pdf_def} leads in fact to the following expression for $L_2$

\begin{equation}
    L_2(p_1,p_2) = d_{1,1} + d_{2,2} - 2\cdot d_{1,2}
    \label{eqn:l2_generic_gmms}
\end{equation}
with

\begin{equation}
    d_{m,n} = \sum_{i=1}^{L_m} \sum_{j=1}^{L_n} \alpha_{m,i}\alpha_{n,j} K(\vb*{\mu}_{m,i},\vb*{\mu}_{n,j},\vb*{P}_{m,i}+\vb*{P}_{n,j}) \qquad m,n = 1,2
\end{equation}
and $K$ defined as

\begin{equation}
    K(\vb*{\mu}_1,\vb*{\mu}_2,\vb*{P}_1+\vb*{P}_2) = \dfrac{1}{\sqrt{|2\pi(\vb*{P}_1+\vb*{P}_2)|}} \cdot \exp{-\frac{1}{2}(\vb*{\mu}_1-\vb*{\mu}_2)^T (\vb*{P}_1+\vb*{P}_2)^{-1} (\vb*{\mu}_1-\vb*{\mu}_2)}
\end{equation}
After setting $p_1=p=p_{\mathcal{N}}(x;0,1)$ and $p_2=\tilde{p}$ equal to the \acrshort{rhs} of \cref{eqn:univariate_gmm_symmetric}, \cref{eqn:l2_generic_gmms} further simplifies into

\begin{equation}
\begin{aligned}
    L_2(p,\tilde{p}) &= 
    \dfrac{1}{2\sqrt{\pi}} - \dfrac{2\tilde{\alpha}_0}{\sqrt{2\pi(1+\tilde{\sigma}^2)}} \\ &+
    \dfrac{1}{2\sqrt{\pi\tilde{\sigma}^2}} \left[\tilde{\alpha}_0^2 + 4\tilde{\alpha}_0\sum_{i=1}^l \tilde{\alpha}_i \exp{-\frac{\tilde{\mu}_i^2}{4\tilde{\sigma}^2}} \right] \\ &+
    \dfrac{1}{\sqrt{\pi\tilde{\sigma}^2}} \sum_{i=1}^l \sum_{j=1}^l \tilde{\alpha}_i \tilde{\alpha}_j \left[\exp{-\frac{(\tilde{\mu}_i-\tilde{\mu}_j)^2}{4\tilde{\sigma}^2}} + \exp{-\frac{(\tilde{\mu}_i+\tilde{\mu}_j)^2}{4\tilde{\sigma}^2}}\right] \\ &-
    \dfrac{4}{\sqrt{2\pi(1+\tilde{\sigma}^2)}} \sum_{i=1}^l \tilde{\alpha}_i \exp{-\frac{\tilde{\mu}_i^2}{2(1+\tilde{\sigma}^2)}}
\end{aligned}
\end{equation}

The nonlinear optimization problem is finally stated as follows.\\
Minimize

\begin{equation}
    J = L_2(p,\tilde{p}) + \lambda\tilde{\sigma}^2
    \label{eqn:split_lib_obj}
\end{equation}
subject to the single equality constraint

\begin{equation}
    0 = 1 - \tilde{\alpha}_0 - 2 \sum_{i=1}^l \tilde{\alpha}_i
    \label{eqn:split_lib_ceq}
\end{equation}
and inequality constraints

\begin{subequations}
    \begin{gather}
        0 < \tilde{\alpha}_i < \tilde{\alpha}_{i-1} \qquad i=2,\ldots,l \\
        0 < \tilde{\mu}_{i-1} < \tilde{\mu}_i \qquad i=2,\ldots,l
    \end{gather}
    \label{eqn:split_lib_neq}
\end{subequations}
where $\lambda > 0$ is a penalty factor that scales the importance of minimizing $\tilde{\sigma}^2$ with respect to $L_2(p,\tilde{p})$. Selecting $\lambda\to 0$ results in less uniform weights $\tilde{\alpha}_i$ with the largest ones allocated to the central kernels, components means $\tilde{\mu}_i$ closer to zero and a larger variance $\tilde{\sigma}^2$. The corresponding splitting libraries will then lead to an higher number of splits when employed by the \acrshort{loads} algorithm. On the contrary, choosing $\lambda\to 1$ results in more equally distributed weights, farther means and a smaller variance, thus reducing the required number of splits.

After selecting a value for $\lambda$, the \acrfull{nlp} resulting from \cref{eqn:split_lib_obj,eqn:split_lib_ceq,eqn:split_lib_neq} is solved numerically using the interior point optimizer \acrshort{ipopt} \cite{Wachter2006} compiled against the \acrshort{hsl} Mathematical Software Library for the solution of sparse linear systems \cite{hsl}. The optimal solution for $L=3$ and $\lambda=10^{-3}$ is given in \cref{tab:split_lib_coef}.

\begin{table}[H]
    \centering
    \begin{tabular}{cccc}
    \toprule
    $i$ & $\tilde{\alpha}_i$ & $\tilde{\mu}_i$ & $\tilde{\sigma}$\\ \midrule
    1 & 0.2252246852539708 & -1.0575150485760967 & 0.6715664864669252 \\
    2 & 0.5495506294920584 & 0.0 & 0.6715664864669252 \\
    3 & 0.2252246852539708 & 1.0575150485760967 & 0.6715664864669252 \\ \bottomrule
    \end{tabular}
    \caption{Univariate splitting library for $L=3$ and $\lambda=10^{-3}$}
    \label{tab:split_lib_coef}
\end{table}

\subsubsection{Splitting Multivariate Distributions}\label{sec:multivariate_splitting}

Gaussian distributions remain Gaussian under linear transformations. Moreover, if $\vb*{X}$ is Gaussian-distributed, then the support of its \acrshort{pdf} is the whole state space $\mathbb{R}^n$. However, for practical purposes, the support can be considered zero outside a \acrshort{ci} around the mean value $\vb*{\mu}_{\vb*{X}}$. For univariate Gaussian distributions, a typical value for its semi-amplitude is that of $3\sigma$ which accounts for $\approx 99.7\%$ of all realizations. Its \acrshort{pdf} may thus be rescaled so as to be zero beyond three standard deviations from its mean. This concept is extended to multivariate \acrshortpl{pdf} by introducing a multivariate splitting procedure that adaptively increases the number of components in \cref{eqn:gmm_pdf_def} and reduces their covariance matrices $\vb*{P}_i$ so that the assumption of a linear transformation holds in the neighborhood of $\vb*{\mu}_i$ in which $p_{\mathcal{N}}(\vb*{x};\vb*{\mu}_i,\vb*{P}_i)$ is considered to be nonzero. This condition must be verified separately for each component $i$ and direction $k$ thus resulting in the recursive algorithm described below.

Given a single Gaussian kernel $p_{\mathcal{N}}(\vb*{x};\vb*{\mu}_i,\vb*{P}_i)$ with weight $\alpha_i$ for which a split is deemed necessary, a univariate splitting library is used to perform this operation in the direction defined by the $k^{th}$ eigenvector of $\vb*{P}_i$ \cite{DeMars2013}. In this work, the splitting direction $k$ is determined by the \acrshort{loads} algorithm described in \cref{sec:loads}. Moreover, the univariate splitting library is computed offline and reused as many time as necessary during the splitting procedure.

Without loss of generality, consider a single kernel of \cref{eqn:gmm_pdf_def} with weight, mean and covariance $\alpha,\vb*{\mu},\vb*{P}$ respectively. The spectral factorization of $\vb*{P}$ is then obtained as

\begin{equation}
    \vb*{P} = \vb*{V}\vb*{\Lambda}\vb*{V}^T
\end{equation}
where $\vb*{V}$ is an orthogonal matrix of right eigenvectors of $\vb*{P}$ and $\vb*{\Lambda}$ a diagonal matrix of corresponding eigenvalues. Performing the split along the $k^{th}$ axis of the spectral factorization results in the following weights, means and covariance for the $L$ new \acrshort{gmm} components

\begin{subequations}
    \begin{align}
        \alpha_i &= \tilde{\alpha_i}\alpha \\
        \vb*{\mu}_i &= \vb*{\mu} + \sqrt{\lambda_k} \tilde{\mu}_i \vb*{v}_k \\
        \vb*{P}_i &= \vb*{V}\vb*{\Lambda}_i\vb*{V}^T
    \end{align}
    \label{eqn:gmm_component_split}
\end{subequations}
with $\vb*{v}_k$ $k^{th}$ column of $\vb*{V}$, $\lambda_k$ corresponding eigenvalue and $\vb*{\Lambda}_i$ defined as

\begin{equation}
    \vb*{\Lambda}_i = \mqty[\dmat{\lambda_1,\ldots,\tilde{\sigma}^2\lambda_k,\ldots,\lambda_{n}}]
    \label{eqn:scaled_eigenvalues_matrix}
\end{equation}

\clearpage
\subsubsection{\glsentrylong{lam} between Distributions}

The \acrfull{lam} is introduced as a measure for the accuracy of the \acrshort{gmm} representation of the transformed \acrshort{pdf}. Given two multivariate \acrshortpl{pdf}, the \acrshort{lam} is defined as \cite{DeMars2013}

\begin{equation}
    \mathcal{L}\left(p,q\right) = \int_{\mathbb{R}^n} p(\vb*{x})q(\vb*{x}) \mathrm{d}\vb*{x}
    \label{eqn:lam_def}
\end{equation}
and describes the amount of overlap between the two distributions, being larger for $p,q$ in greater agreement between each other.

If the agreement of a \acrshort{gmm} with respect to samples resulting from a \acrshort{mc} simulation is of interest, \cref{eqn:lam_def} can be reformulated treating the \acrshort{mc} samples as a \acrfull{dmm} as

\begin{equation}
    \mathcal{L}\left(p,q\right) = \sum_{i=1}^{N_s}\sum_{j=1}^{L} \gamma_i\alpha_j p_{\mathcal{N}}(\vb*{x}_i;\vb*{\mu}_j,\vb*{P}_j)
    \label{eqn:lam_dmm}
\end{equation}
with $N_s,L$ number of \acrshort{mc} samples and Gaussian kernels respectively. Moreover, $\alpha_j,\vb*{\mu}_j,\vb*{P}_j$ are the \acrshort{gmm} components weights, means and covariances while $\vb*{x}_i,\gamma_i=1/N_s\ \forall i$ are the propagated random samples and corresponding weights.\par

\subsection{Combined \glsentryshort{loads}-\glsentryshort{gmm} algorithm}\label{sec:loads_gmm_algorithm}

In this work, a univariate splitting library such as the one in \cref{tab:split_lib_coef} is embedded within the \acrshort{loads} framework and the \acrshortpl{nli} $\nu,\nu_e$ defined in \cref{eqn:nonlinear_index_def,eqn:directional_nonlinear_index_def} are used to identify the need for new splits and the corresponding directions. Moreover, an additional constraint is taken into account to bound the components weights, i.e. $\alpha^{(p)} \geq \alpha_{min}\ \forall p$. Given the nonlinear transformation in \cref{eqn:generic_nonlinear_transform}, the developed algorithm is as follows.\par

Firstly, the initial uncertainty on the function input $\vb*{x}$ is modeled with a multivariate Gaussian distribution with known mean $\vb*{\mu}_0$ and covariance $\vb*{P}_0$. A \acrshort{ci} $c$ is then chosen and the input \acrshort{da} state is initialized as

\begin{equation}
    [\vb*{x}] = \vb*{\mu}_0 + \vb*{V}_0 \cdot \left\{c\sqrt{\vb*{\lambda}_0}\cdot\delta\vb*{x}\right\} = \bar{\vb*{x}}_0 + \vb*{\beta}_0\cdot\delta\vb*{x}
    \label{eqn:da_input}
\end{equation}
where $\vb*{\lambda}_0,\vb*{V}_0$ are the eigenvalues and eigenvectors of $\vb*{P}_0$ such that $\vb*{P}_0 = \vb*{V}_0\vb*{\Lambda}_0\vb*{V}_0^T$ with $\vb*{\Lambda}_0=diag(\vb*{\lambda}_0)$. A single \acrshort{gmm} component with $\alpha=1,\vb*{\mu}_0,\vb*{P}_0$ is then associated to $[\vb*{x}]$.

Secondly, a univariate splitting library $\mathcal{S}_l$ is selected and the steps of the iterative \acrshort{loads} algorithm described in \cref{sec:loads} reformulated as in \cref{alg:loads_gmm_algorithm}.

\begin{algorithm}[H]
\begin{algorithmic}
\setstretch{1.2}
    \Function{loadsGmm}{$\vb*{f},[\vb*{x}],\alpha_0,\vb*{\mu}_0,\vb*{P}_0,\alpha_{min},\varepsilon_\nu,N_{max},\mathcal{S}_l$}
        \State Initialize the working and output manifolds $W_{\vb*{x}},M_{\vb*{x}},M_{\vb*{y}}$
        \State Add $[\vb*{x}]$ to $W_{\vb*{x}}$
        \While{$W_{\vb*{x}}$ is not empty}
            \State Remove the first set $[\vb*{x}^{(p)}]$ from $W_{\vb*{x}}$
            \State Evaluate $[\vb*{y}^{(p)}] = \vb*{f}([\vb*{x}^{(p)}])$
            \If{$\alpha^{(p)} < \alpha_{min}$}
                \State Add $[\vb*{x}^{(p)}]$ to $M_{\vb*{x}}$ and $[\vb*{y}^{(p)}]$ to $M_{\vb*{y}}$
            \Else
                \State Compute the scaled Jacobian $[\vb*{J}^{(p)}]$
                \begin{equation}
                    [J_{ij}^{(p)}] = \dfrac{1}{c\sqrt{\lambda_j^{(p)}}}\cdot\dfrac{\partial f_i}{\partial x_j} \qquad i=1,\ldots,m \quad j=1,\ldots,n
                    \label{eqn:jacobian_step}
                \end{equation}
                \State Compute the \acrshort{nli} $\nu^{(p)}$ \Comment{see \cref{eqn:jacobian_component,eqn:nonconstant_jacobian_component,eqn:jacobian_component_bounds,eqn:nonlinear_index_def,eqn:jacobian_step}}

                \If{$\nu^{(p)} < \varepsilon_{\nu} \lor N_{splits}^{(p)}=N_{max}$}
                    \State Add $[\vb*{x}^{(p)}]$ to $M_{\vb*{x}}$ and $[\vb*{y}^{(p)}]$ to $M_{\vb*{y}}$
                \Else
                    \State Compute the splitting direction $k$ \Comment{see \cref{eqn:directional_jacobian,eqn:directional_nonlinear_index_def,eqn:splitting_direction,eqn:jacobian_step}}
                    \State Split $[\vb*{x}^{(p)}]$ into $[\vb*{x}^{(i)}]$ for $i=1,\ldots,L$
                    \begin{equation}
                        [\vb*{x}^{(i)}] = [\vb*{x}^{(p)}] \circ \left\{0,\ldots,\dfrac{\tilde{\mu}_i}{L} + \tilde{\sigma}\delta x_k,\ldots,0\right\}^T
                    \end{equation}
                    \State Split $(\alpha^{(p)},\vb*{\mu}^{(p)},\vb*{P}^{(p)})$ into $(\alpha^{(i)},\vb*{\mu}^{(i)},\vb*{P}^{(i)})$ for $i=1,\ldots,L$ \Comment{see \cref{eqn:gmm_component_split}}
                    \begin{subequations}
                        \begin{align}
                            \alpha^{(i)} &= \tilde{\alpha_i}\alpha^{(p)} \\
                            \vb*{\mu}^{(i)} &= \vb*{\mu}^{(p)} + \sqrt{\lambda_k^{(p)}} \tilde{\mu}_i \vb*{v}_k^{(p)} \\
                            \vb*{P}^{(i)} &= \vb*{V}^{(p)}\vb*{\Lambda}_i^{(p)}\vb*{V}^{(p),T}
                        \end{align}
                    \end{subequations}
                    \State Add all $[\vb*{x}^{(i)}]$ and $(\alpha^{(i)},\vb*{\mu}^{(i)},\vb*{P}^{(i)})$ to $W_{\vb*{x}}$
                \EndIf
            \EndIf
        \EndWhile
        \State \Return $M_{\vb*{x}},M_{\vb*{y}}$
    \EndFunction
\end{algorithmic}
\caption{Combined \glsentryshort{loads}-\glsentryshort{gmm} algorithm}
\label{alg:loads_gmm_algorithm}
\end{algorithm}

Once \cref{alg:loads_gmm_algorithm} has converged, the \acrshort{pdf} of the initial random variable $\vb*{X}\sim\mathcal{N}(\vb*{\mu}_0,\vb*{P}_0)$ is approximated by a multivariate \acrshort{gmm} whose components weights, means and covariances are stored in $M_{\vb*{x}}$, namely

\begin{equation}
    p_{\vb*{X}}(\vb*{x}) \simeq \sum_{p=1}^{P} \alpha^{(p)} p_{\mathcal{N}}\left(\vb*{x};\vb*{\mu}^{(p)},\vb*{P}^{(p)}\right)
    \label{eqn:split_initial_gmm}
\end{equation}
with $P$ dimension of $M_{\vb*{x}},M_{\vb*{y}}$. One more step is thus needed to estimate $p_{\vb*{Y}}$ as described in \cref{sec:gmm_pdf_propagation}.

\subsubsection{Estimating \glsentryshort{gmm} \glsentryshortpl{pdf}}\label{sec:gmm_pdf_propagation}

Outputs of the \acrshort{loads}-\acrshort{gmm} algorithm described above comprise the \acrshort{gmm} approximation of the initial \acrshort{pdf} given by \cref{eqn:split_initial_gmm} and the Taylor expansions of the function image $[\vb*{y}^{(p)}] = \vb*{f}([\vb*{x}^{(p)}])$ stored in $M_{\vb*{y}}$. The transformed \acrshort{pdf} $p_{\vb*{Y}}$ is then estimated as follows.

Firstly, a sampling scheme is chosen to draw samples from each component of the initial distribution in \cref{eqn:split_initial_gmm}. Several sampling schemes such as \acrfull{mc}, \acrfull{ut} \cite{Julier2004} and \acrfull{cut} \cite{Adurthi2018} might be employed. In this work, the \acrshort{ut} scheme is selected since it minimizes the number sigma points while guaranteeing an exact representation of the first two statistical moments of the reconstructed \acrshort{pdf}. This choice is deemed consistent with the second-order Taylor expansions available for $[\vb*{y}^{(p)}]$. Following similar arguments, \acrshort{cut4} and \acrshort{cut6} schemes should be preferred when higher-order expansions are instead available in $M_{\vb*{y}}$. It was also verified that \acrshort{mc} sampling does not provide any significant improvement in the accuracy of the estimated \acrshort{pdf} with respect to the selected \acrshort{ut} scheme.

Secondly, each \acrshort{gmm} component is sampled according to its mean $\vb*{\mu}^{(p)}$ and covariance $\vb*{P}^{(p)}$. The drawn samples are then mapped through $\vb*{f}$ evaluating the corresponding Taylor expansion of the output $[\vb*{y}^{(p)}] = \vb*{f}([\vb*{x}^{(p)}])$ stored in $M_{\vb*{y}}$. The first two moments of the transformed \acrshort{pdf} are then obtained from the transformed samples as

\begin{subequations}
    \begin{align}
        \vb*{\mu}_{\vb*{Y}}^{(p)} &= \frac{1}{N_s} \sum_{s=1}^{N_s} \vb*{y}_s^{(p)}\\
        \vb*{P}_{\vb*{Y}}^{(p)} &= \frac{1}{N_s-1} \sum_{s=1}^{N_s} \qty(\vb*{y}_{s}^{(p)} - \vb*{\mu}^{(p)})\qty(\vb*{y}_{s}^{(p)} - \vb*{\mu}^{(p)})^T
    \end{align}
    \label{eqn:moments_mc}
\end{subequations}
for \acrshort{mc} samples and

\begin{subequations}
    \begin{align}
        \vb*{\mu}_{\vb*{Y}}^{(p)} &= \sum_{s=1}^{N_s} w_s \vb*{y}_s^{(p)}\\
        \vb*{P}_{\vb*{Y}}^{(p)} &= \sum_{s=1}^{N_s} w_s \qty(\vb*{y}_s^{(p)} - \vb*{\mu}^{(p)}) \qty(\vb*{y}_s^{(p)} - \vb*{\mu}^{(p)})^T
    \end{align}
    \label{eqn:moments_ut}
\end{subequations}
for \acrshort{ut} and \acrshort{cut} points. Here $N_s$ is the total number of samples and $w_s$ the associated weights with $\sum_{s=1}^{N_s} w_s = 1$.\par

The transformed \acrshort{pdf} is finally estimated as

\begin{equation}
    p_{\vb*{Y}}(\vb*{y}) \simeq \sum_{p=1}^{P} \alpha^{(p)} p_{\mathcal{N}}\left(\vb*{y};\vb*{\mu}^{(p)},\vb*{P}^{(p)}\right)
    \label{eqn:transformed_pdf_gmm}
\end{equation}
with $\vb*{\mu}_{\vb*{Y}}^{(p)},\vb*{P}_{\vb*{Y}}^{(p)}$ given by \cref{eqn:moments_mc} or \cref{eqn:moments_ut}.

\section{\glsentrytitlecase{mf}{long} Orbit \glsentrylong{up}}\label{sec:mf_orbit_up}

The combined \acrshort{loads}-\acrshort{gmm} algorithm presented in \cref{sec:loads_gmm_algorithm} can be applied to any nonlinear transformation $\vb*{f}$. In the followings, its application within a \acrfull{mf} method for \acrshort{up} in orbital dynamics is presented. The different \acrfullpl{oe} used to describe the orbit state are firstly introduced in \cref{sec:orbital_elements} followed by a description of the \acrshort{da}-based \acrfull{sgp} dynamical model and the formulation of the \acrlong{mf} \acrshort{loads}-\acrshort{gmm} algorithm in \cref{sec:da_sgp4,sec:mf_loads_gmm} respectively.\par

\subsection{\glsentrylongpl{oe}}\label{sec:orbital_elements}

Four \acrfullpl{oe} sets are introduced in this section. Following the work of \citet{Junkins2004}, their impact on the overall performance of the proposed \acrlong{mf} \acrshort{up} method is assessed in \cref{sec:numerical_results}.

\subsubsection{Cartesian Parameters}\label{sec:cartesian}

Cartesian parameters describe the spacecraft position and velocity in three-dimensional Euclidean space, i.e. $\vb*{x} = [\vb*{r}\ \vb*{\dot{r}}]^T = [x\ y\ z\ \dot{x}\ \dot{y}\ \dot{z}]^T$. These parameters are well suited to describe third-body perturbations and nonconservative forces (\acrshort{srp}, atmospheric drag) acting on the spacecraft and they are commonly used in space surveillance for collision risk assessment. Under perturbed Keplerian dynamics, their evolution over time is governed by the second-order \acrshort{ode}

\begin{equation}
    \vb*{\ddot{r}} = -\dfrac{\mu}{r^3}\vb*{r} + \vb*{a}_p(\vb*{r},\vb*{\dot{r}},t)
    \label{eqn:cartesian_odes}
\end{equation}

with $\mu$ central body standard gravitational parameter, $r=\norm{\vb*{r}}_2$ and $\vb*{a}_p$ sum of all perturbing accelerations other than two-body central attraction. \cref{eqn:cartesian_odes} is strongly nonlinear in the cartesian parameters due to the Keplerian attraction being proportional to $1/r^2$. Departure from Gaussianity will thus be magnified by the choice of this particular coordinates set resulting in a larger number of Gaussian kernels required to capture the true state \acrshort{pdf}.

\subsubsection{Equinoctial \glsentrylongpl{oe}}\label{sec:equinoctial}

Equinoctial \acrshortpl{oe} were introduced by \citet{Broucke1972} as an alternative to the well known Keplerian parameters to overcome some issues intrinsic of the classical set, namely a singularity for circular and/or equatorial orbits for which the \acrlong{aop} and/or the line of nodes is not defined. These new elements are expressed in terms of the Keplerian ones as\footnote{the notation and order introduced in \citet{Walker1985} is used here instead of that in \cite{Broucke1972}.}

\begin{equation}
    \begin{aligned}
        a &= a\\
        f &= e \cos{\left(\omega + \Omega\right)}\\
        g &= e \sin{\left(\omega + \Omega\right)}\\
        h &= \tan{\left(\frac{i}{2}\right)} \cos{\left(\Omega\right)}\\
        k &= \tan{\left(\frac{i}{2}\right)} \sin{\left(\Omega\right)}\\
        \lambda_0 &= \Omega + \omega + M_0
    \end{aligned}
    \label{eqn:equinoctial_elements_def}
\end{equation}

with $a$ semimajor axis, $e$ eccentricity, $i$ inclination, $\Omega$ \acrfull{raan}, $\omega$ \acrfull{aop} and $M_0$ mean anomaly. In this work, the mean longitude $\lambda$ or true longitude $L$ are used instead of $\lambda_0$. They are defined as

\begin{equation}
    \begin{aligned}
        \lambda &= \Omega + \omega + M\\
        L &= \Omega + \omega + \nu
    \end{aligned}
    \label{eqn:mean_true_longitudes_def}
\end{equation}

with $M,\nu$ mean and true anomalies respectively\footnote{$\lambda,\nu$ here denote the mean longitude and true anomaly not to be confused with the matrix eigenvalues and the \acrlong{nli} of \cref{eqn:nonlinear_index_def}.}.

The set of \acrshortpl{oe} in \cref{eqn:equinoctial_elements_def} is nonsingular except for retrograde equatorial orbits $(i=\pi)$ and is well suited for the numerical propagation of perturbed Keplerian dynamics. Yet its application is limited to elliptic orbits $(e<1)$ with positive semimajor axis $(a>0)$. Gauss planetary equations are used throughout this paper to numerically propagate orbit states expressed in equinoctial \acrshortpl{oe} \cite{Battin1999}.

\subsubsection{\glsentrylongpl{mee}}\label{sec:modified_equinoctial}

\acrfullpl{mee} were firstly proposed by \citet{Walker1985} as an \acrshortpl{oe} set suitable for perturbation analysis of all kind of orbits. The only difference with respect to equinoctial \acrshortpl{oe} is the substitution of the semimajor axis $a$ by the semilatus rectum $p=a(1-e^2)$ and the mean longitude at epoch $\lambda_0$ by the true longitude $L$ given by \cref{eqn:mean_true_longitudes_def}. The remaining four elements are $f,g,h,k$ defined in \cref{eqn:equinoctial_elements_def}.

\acrshortpl{mee} are applicable to all type of orbits (elliptic, parabolic, hyperbolic) excluding retrograde equatorial ones $(i=\pi)$ for which they also exhibit a singularity. These elements are particularly suitable to describe optimal control problems arising in low-thrust transfers design. As for the equinoctial \acrshortpl{oe}, Gauss planetary equations are employed in this work to describe the time evolution of \acrshortpl{mee} while the mean longitude $\lambda$ is used interchangeably with $L$.

\subsubsection{Alternate Equinoctial Elements}\label{sec:alternate_equinoctial}

Alternate equinoctial elements were derived from the equinoctial \acrshortpl{oe} by \citet{Horwood2011} substituting $a,\lambda_0$ by $n,\lambda$ where $n=\sqrt{\mu/a^3}$ is the Keplerian mean motion and $\lambda$ the mean longitude defined in \cref{eqn:mean_true_longitudes_def}. The other four elements are $f,g,h,k$ as in \cref{eqn:equinoctial_elements_def}\footnote{the notation and order introduced in \cite{Walker1985} is used here instead of that in \cite{Horwood2011}.}. Alternate equinoctial elements evolve linearly in time under the assumption of unperturbed Keplerian motion ($n,f,g,h,k$ constant and $d\lambda/dt=n$) and nearly linearly in presence of perturbing forces. In the context of orbit \acrshort{up}, these elements will therefore preserve Gaussianity for longer propagation time spans compared to other coordinates sets. In this work, modified Gauss planetary equations that uses $n$ as first coordinate instead of either $a$ or $p$ are employed to describe their evolution over time.

\subsection{\glsentryshort{da} based \glsentrylong{sgp} Model}\label{sec:da_sgp4}

The \acrfull{sgp} model series were originally developed by the \acrshort{us} \acrfull{dod} for space surveillance applications \cite{Vallado2006a}. Five models have been released over the years, namely \acrshort{sgp}, \acrshort{sgp}4, \acrshort{sdp}4, \acrshort{sgp}8 and \acrshort{sdp}8. Nowadays, \acrshort{sgp}4/\acrshort{sdp}4 are collectively known as \acrshort{sgp}4 theory and commonly used for the propagation of \acrfullpl{tle}. \acrshort{sgp}4 is used for near-Earth objects (orbital period $<225\ min$) and include the effects of Earth's nonuniform gravity field and atmospheric drag. \acrshort{sdp}4 is used for deep-space objects (orbital period $\geq225\ min$) and adds to \acrshort{sgp}4 the contribution of lunisolar perturbations and Earth's resonance terms for geosynchronous and highly-elliptical semi-synchronous orbits. \acrshortpl{tle} produced daily by \acrshort{norad} contain the mean \acrshortpl{oe} of a space object and other relevant information. These elements are intrinsic to the \acrshort{sgp}4 theory and osculating \acrshortpl{oe} at any epoch must be computed with the relevant \acrshort{sgp}4/\acrshort{sdp}4 propagator. The lasts are analytical functions that computes an osculating state in cartesian parameters given a \acrshort{tle} and the epoch $\tau$ at which the state is sought

\begin{equation}
    \vb*{x}_{TEME}(\tau) = \vb*{f}_{SGP4}(\tau,TLE)
    \label{eqn:sgp4_evaluation}
\end{equation}

Mean \acrshort{tle} elements and osculating states output of \cref{eqn:sgp4_evaluation} are expressed in \acrfull{teme} frame. Since an exact definition of this frame is not available in the literature, it is recommended to convert \acrshort{teme} states into a truly standard coordinate system before interfacing with external routines \cite{Vallado2006a}. Transformations from \acrshort{teme} to \acrfull{gcrf} are described for example in \citet{Vallado2013}.

If the computation of $\vb*{x}_{TEME}$ from a \acrshort{tle} and $\tau$ is straightforward and consists in a single function evaluation, the opposite is not true and care must be taken while converting an osculating state into the equivalent \acrshort{tle} representation. Inversion of \cref{eqn:sgp4_evaluation} can still be accomplished with the iterative procedure summarized in \cref{alg:state2tle}.

\begin{algorithm}[ht]
\begin{algorithmic}
    \setstretch{1.3}
    \Function{stateToTLE}{$\tau,\vb*{x}_{TEME},\varepsilon_{max},i_{max}$}
        \State Compute osculating Keplerian elements $\vb*{x}_{kep,osc}$ from $\vb*{x}_{TEME}$
        \State Initialize $\vb*{x}_{cart}^{(i)}$ with $\vb*{x}_{TEME}$
        \State Initialize $\vb*{x}_{kep,osc}^{(i)}$ with $\vb*{x}_{kep,osc}$
        \State Initialize mean elements $\vb*{x}_{kep,mean}^{(i)}$ with $\vb*{x}_{kep,osc}^{(i)}$
        \State $i \gets 0$
        \While {$i < i_{max}$}
            \State Evaluate \cref{eqn:sgp4_evaluation} in $\tau,\vb*{x}_{kep,mean}^{(i)}$ and update $\vb*{x}_{cart}^{(i)}$
            \State Compute $\vb*{x}_{kep,osc}^{(i)}$ from $\vb*{x}_{cart}^{(i)}$
            \State Compute $\Delta\vb*{x}_{kep,osc}^{(i)}=\vb*{x}_{kep,osc}^{(i)}-\vb*{x}_{kep,osc}$
            \If {$\norm*{\Delta\vb*{x}_{kep,osc}^{(i)}}<\varepsilon_{max}$}
                \State Build the \acrshort{tle} from $\vb*{x}_{kep,mean}^{(i)}$
                \State \textbf{break}
            \Else
                \State Compute $\vb*{x}_{kep,mean}^{(i)}=\vb*{x}_{kep,mean}^{(i)}+\Delta\vb*{x}_{kep,osc}^{(i)}$
            \EndIf
            \State $i \gets i+1$
        \EndWhile
        \State \Return TLE
        \EndFunction
\end{algorithmic}
\caption{\glsentryshort{tle} from osculating cartesian state}
\label{alg:state2tle}
\end{algorithm}

Evaluation of \cref{eqn:sgp4_evaluation} and \cref{alg:state2tle} do not pose any relevant computational issue if performed in \acrfull{fp} arithmetic. Moreover, \cref{eqn:sgp4_evaluation} contains only algebraic operations that are easily implemented in the \acrshort{da} framework. However, \cref{alg:state2tle} cannot be directly applied to Taylor polynomials due to the iterative procedure and associated convergence criterion. In this work, a simple check on the constant part of $[\Delta\vb*{x}_{kep,osc}^{(i)}]$ is performed to exit the while loop in \cref{alg:state2tle}. $\norm*{\Delta\vb*{x}_{kep,osc}^{(i)}}$ is thus replaced by $\norm*{\overline{\Delta\vb*{x}}_{kep,osc}^{(i)}}=\norm*{cons([\Delta\vb*{x}_{kep,osc}^{(i)}])}$.

Denoting with $[\vb*{x}_{TEME}] = \overline{\vb*{x}}_{TEME} + \vb*{\beta}\cdot\delta\vb*{x}_{TEME}$ the nominal cartesian state and associated uncertainty, the evaluation of \cref{alg:state2tle} on $(\tau,[\vb*{x}_{TEME}])$ outputs the Taylor expansion of the \acrshort{tle} with respect to deviations in the initial osculating state, i.e. $[TLE]=\mathcal{T}_{TLE}(\delta\vb*{x}_{TEME})$. The uncertainty in the mean \acrshort{tle} elements is thus obtained directly from $\mathcal{T}_{TLE}(\delta\vb*{x}_{TEME})$ while a posteriori estimation is required if \acrshort{fp} arithmetic is used.

Choosing \acrshort{sgp}4 as \acrlong{lf} dynamical model for the \acrlong{mf} \acrshort{up} method ultimately relies on its integration within the \acrshort{loads}-\acrshort{gmm} framework described in \cref{sec:uncertainty_propagation}. A function in the form of \cref{eqn:generic_nonlinear_transform} is thus sought in which both $\vb*{x}$ and $\vb*{y}$ are osculating orbit states expressed in either one of the \acrshortpl{oe} sets in \cref{sec:orbital_elements}. A wrapper is then built around \cref{eqn:sgp4_evaluation} to be used as \acrshort{loads} target function as summarized in \cref{alg:sgp4_target}.

\begin{algorithm}[H]
\begin{algorithmic}
    \setstretch{1.3}
    \Function{targetSGP4}{$\tau_0,\tau,[\vb*{x}_{osc}(\tau_0)]$}
        \State Compute cartesian state $[\vb*{x}_{cart}(\tau_0)]$ from $[\vb*{x}_{osc}(\tau_0)]$
        \State Transform $[\vb*{x}_{cart}(\tau_0)]$ into $[\vb*{x}_{TEME}(\tau_0)]$ expressed in \acrshort{teme} frame
        \State Compute initial \acrshort{tle} from $[\vb*{x}_{TEME}(\tau_0)]$ \Comment{see \cref{alg:state2tle}}
        \State Compute $[\vb*{x}_{TEME}(\tau)]$ from initial TLE \Comment{see \cref{eqn:sgp4_evaluation}}
        \State Transform $[\vb*{x}_{TEME}(\tau)]$ into $[\vb*{x}_{cart}(\tau)]$ expressed in original frame
        \State Compute final state $[\vb*{x}_{osc}(\tau)]$ from $[\vb*{x}_{cart}(\tau)]$
        \State \Return $[\vb*{x}_{osc}(\tau)]$
    \EndFunction
\end{algorithmic}
\caption{\glsentryshort{loads} target function for \glsentryshort{sgp}4 dynamical model}
\label{alg:sgp4_target}
\end{algorithm}

\subsection{\glsentrytitlecase{mf}{long} \glsentryshort{loads}-\glsentryshort{gmm} Algorithm}\label{sec:mf_loads_gmm}

The \acrshort{mf} method developed in this work comprises two steps, a \acrfull{lf} propagation carried out in \acrshort{da} framework and an \acrfull{hf} one executed \acrfull{pw} for a posteriori correction of the \acrshort{lf} step.

The \acrshort{sgp}4 model described in \cref{sec:da_sgp4} is chosen for the \acrshort{lf} dynamics since its analytical formulation allows a straightforward integration with the \acrshort{loads}-\acrshort{gmm} algorithm. An accurate numerical propagator is then used for the \acrshort{hf} correction. The overall \acrshort{mf} method is detailed in \cref{alg:mf_orbit_up}.

\begin{algorithm}[H]
\begin{algorithmic}
\setstretch{1.3}
    \Function{multifidelityOrbitUP}{$t_0,t,\vb*{\mu}_{\vb*{X}}(t_0),\vb*{P}_{\vb*{X}}(t_0),c,\alpha_{min},\varepsilon_\nu,N_{max},\mathcal{S}_l$}
    
        \State Model the initial uncertainty on the orbit state as a multivariate Gaussian distribution
        \begin{equation}
            p_{\vb*{X}}(\vb*{x}(t_0))\sim\mathcal{N}\left(\vb*{\mu}_{\vb*{X}}(t_0),\vb*{P}_{\vb*{X}}(t_0)\right)
            \label{eqn:initial_normal_distribution_mf}
        \end{equation}
        
        \State Initialize $[\vb*{x}]=[\vb*{x}_{LF}(t_0)]$ from $\vb*{\mu}_{\vb*{X}}(t_0),\vb*{P}_{\vb*{X}}(t_0)$ \Comment{see \cref{eqn:da_input}}
        \State Run the \acrshort{loads}-\acrshort{gmm} algorithm with \acrshort{sgp}4 as target function $\vb*{f}$ \Comment{see \cref{alg:loads_gmm_algorithm,alg:sgp4_target}}

        \State Approximate the initial \acrshort{pdf} $p_{\vb*{X}}(\vb*{x}(t_0))$ \Comment{see \cref{eqn:split_initial_gmm}}
        \begin{equation}
            p_{\vb*{X}}(\vb*{x}(t_0)) \simeq \sum_{p=1}^{P} \alpha^{(p)} p_{\mathcal{N}}\left(\vb*{x};\vb*{\mu}_{LF}^{(p)}(t_0),\vb*{P}_{LF}^{(p)}(t_0)\right)
        \end{equation}
        
        \State Propagate the initial means $\vb*{\mu}_{LF}^{(p)}(t_0)$ in \acrlong{hf} to obtain $\vb*{\mu}_{HF}^{(p)}(t)$
        
        \State Separate the constant part and nilpotent part of the Taylor expansions $[\vb*{x}_{LF}^{(p)}(t)]$
        \begin{equation}
            [\vb*{x}_{LF}^{(p)}(t)] = \overline{\vb*{x}}_{LF}^{(p)}(t) + \mathcal{T}_{\delta\vb*{x}_{LF}^{(p)}(t)}(t_0,\delta\vb*{x}_0)
            \label{eqn:lf_expansion_splitted}
        \end{equation}
        
        \State Substitute $\overline{\vb*{x}}_{LF}^{(p)}(t)$ in \cref{eqn:lf_expansion_splitted} with $\vb*{\mu}_{HF}^{(p)}(t)$
        \begin{equation}
            [\vb*{x}_{MF}^{(p)}(t)] = \vb*{\mu}_{HF}^{(p)}(t) + \mathcal{T}_{\delta\vb*{x}_{LF}^{(p)}(t)}(t_0,\delta\vb*{x}_0)
            \label{eqn:mf_expansion_splitted}
        \end{equation}
        
        \State Estimate the transformed \acrshort{pdf} from $[\vb*{x}_{MF}^{(p)}(t)]$ \Comment{see \cref{eqn:split_initial_gmm,eqn:transformed_pdf_gmm}}
        \begin{equation}
            p_{\vb*{X}}(\vb*{x}(t)) \simeq \sum_{p=1}^{P} \alpha^{(p)} p_{\mathcal{N}}\left(\vb*{x};\vb*{\mu}_{MF}^{(p)}(t),\vb*{P}_{MF}^{(p)}(t)\right)
        \end{equation}
        
        \State \Return $p_{\vb*{X}}(\vb*{x}(t))$
    \EndFunction
\end{algorithmic}
\caption{\glsentrytitlecase{mf}{long} orbit \glsentrylong{up} algorithm}
\label{alg:mf_orbit_up}
\end{algorithm}

If \acrshort{sgp}4 is used as \acrshort{lf} dynamical model, the polynomial shift in \cref{eqn:mf_expansion_splitted} might be performed in either osculating or mean \acrshort{tle} elements space. Being $\vb*{\mu}_{HF}^{(p)}(t)$ and $\mathcal{T}_{\delta\vb*{x}_{LF}^{(p)}(t)}(t_0,\delta\vb*{x}_0)$ expressed in the chosen osculating coordinates set, a single step is required to apply \cref{eqn:mf_expansion_splitted} in this space. Conversely, intermediate operations are needed to perform this shift in \acrshort{tle} space. The procedure is detailed in \cref{alg:mf_shift_tle} while its effect on the overall accuracy of the \acrshort{up} method is discussed in \cref{sec:numerical_results}.

\begin{algorithm}[H]
\begin{algorithmic}
\setstretch{1.3}
    \Function{shiftTLE}{$\vb*{\mu}_{HF}^{(p)}(t),\mathcal{T}_{\delta\vb*{x}_{LF}^{(p)}(t)}(t_0,\delta\vb*{x}_0)$}
    
        \State Transform $\vb*{\mu}_{HF}^{(p)}(t),\mathcal{T}_{\delta\vb*{x}_{LF}^{(p)}(t)}(t_0,\delta\vb*{x}_0)$ into $\vb*{\mu}_{TEME,HF}^{(p)}(t),\mathcal{T}_{\delta\vb*{x}_{TEME,LF}^{(p)}(t)}(t_0,\delta\vb*{x}_0)$
        
        \State Compute $\vb*{\mu}_{TLE,HF}^{(p)}(t),\mathcal{T}_{\delta\vb*{x}_{TLE,LF}^{(p)}(t)}(t_0,\delta\vb*{x}_0)$ from $\vb*{\mu}_{TEME,HF}^{(p)}(t),\mathcal{T}_{\delta\vb*{x}_{TEME,LF}^{(p)}(t)}(t_0,\delta\vb*{x}_0)$ \Comment{see \cref{alg:state2tle}}
    
        \State Compute $[\vb*{x}_{TLE,MF}^{(p)}(t)] = \vb*{\mu}_{TLE,HF}^{(p)}(t) + \mathcal{T}_{\delta\vb*{x}_{TLE,LF}^{(p)}(t)}(t_0,\delta\vb*{x}_0)$
        
        \State Compute $[\vb*{x}_{TEME,MF}^{(p)}(t)]$ from $[\vb*{x}_{TLE,MF}^{(p)}(t)]$ \Comment{see \cref{eqn:sgp4_evaluation}}
        
        \State Transform $[\vb*{x}_{TEME,MF}^{(p)}(t)]$ into $[\vb*{x}_{MF}^{(p)}(t)]$
        
        \State \Return $[\vb*{x}_{MF}^{(p)}(t)]$
    \EndFunction
\end{algorithmic}
\caption{Polynomial shift in mean \glsentryshort{tle} elements space}
\label{alg:mf_shift_tle}
\end{algorithm}

\section{Numerical Results}\label{sec:numerical_results}

Three applications of the \acrlong{mf} method for orbit \acrlong{up} are presented in this section. The analytical \acrshort{sgp}4 model is used for the \acrshort{lf} dynamics while the \acrshort{hf} propagation is carried out numerically taking into account the following force models

\begin{enumerate}
    \item Earth non-uniform gravity field with $8 \times 8$ spherical harmonics potential \cite{Holmes2002}
    \item Sun and Moon third-body forces with bodies positions from \acrshort{jpl} ephemerides
    \item Isotropic drag force with modified Harris-Priester static atmosphere model \cite{Harris1962} (\acrshort{leo} case only)
    \item Isotropic \acrfull{srp} force including Earth's umbra and penumbra transitions \cite{Hubaux2012}
\end{enumerate}

The Dormand-Prince 8(5,3) integrator available in Hipparchus\footnote{\href{https://hipparchus.org/}{https://hipparchus.org/}} is used for the numerical propagation of the \acrlong{hf} \acrshortpl{ode} while the Earth potential and atmosphere models are those of Orekit\footnote{\href{https://www.orekit.org/}{https://www.orekit.org/}}. All simulations were run with initial epoch set on January 1\textsuperscript{st}, 2021 00:00:00 \acrshort{utc} and force models parameters summarized in \cref{tab:forces_params}.

\begin{table}[H]
    \setstretch{1.2}
    \centering
    \begin{tabular}{ccl}
        \toprule
        Parameter & Value & Description\\ \midrule
        $\mu_{E}$ & \SI{398600.4355}{\cubic\km\per\square\s} & Earth standard gravitational parameter\\
        $m$ & \SI{500.0}{\kg} & satellite mass\\
        $A_n$ & \SI{1.0}{\square\m} & satellite cross section area\\
        $C_D$ & 2.0 & atmospheric drag coefficient\\
        $C_R$ & 1.5 & \acrshort{srp} reflection coefficient\\
        \bottomrule
    \end{tabular}
    \caption{Force models parameters}
    \label{tab:forces_params}
\end{table}

The first two simulations are adapted from \cite{DeMars2013} and include an \acrfull{heo} and a \acrfull{leo} cases. The last scenario is a \acrfull{meo} orbit for the Galileo constellation. Nominal \acrfullpl{ic} are summarized in \cref{tab:nominal_ics}. The corresponding initial uncertainties are assumed Gaussian with $1\sigma$ deviations in cartesian parameters given in \cref{tab:initial_sigma}.

\begin{table}[H]
    \setstretch{1.2}
    \centering
    \begin{tabular}{ccccccc}
    \toprule
    Regime & $a$, \si{\km} & $e$ & $i$, \si{\deg} & $\Omega$, \si{\deg} & $\omega$, \si{\deg} & $M$, \si{\deg} \\ \midrule
    \acrshort{heo} & \num{35000.0} & 0.2 & 0.0 & 0.0 & 0.0 & 0.0 \\
    \acrshort{leo} & \num{6678.0} & 0.01 & 0.0 & 0.0 & 0.0 & 0.0 \\
    \acrshort{meo} & \num{29600.135} & 0.0 & \num{56.0} & 0.0 & 0.0 & 0.0 \\
    \bottomrule
    \end{tabular}
    \caption{Nominal \glsentrylongpl{ic} for all test cases}
    \label{tab:nominal_ics}
\end{table}

\begin{table}[H]
    \setstretch{1.2}
    \centering
    \begin{tabular}{ccccccc}
    \toprule
    Regime & $\sigma_x$, \si{\km} & $\sigma_y$, \si{\km} & $\sigma_z$, \si{\km} & $\sigma_{vx}$, \si{\m\per\s} & $\sigma_{vy}$, \si{\m\per\s} & $\sigma_{vz}$, \si{\m\per\s} \\ \midrule
    \acrshort{heo} & 1.0 & 1.0 & 0.0 & 1.0 & 1.0 & 0.0 \\
    \acrshort{leo} & 1.3 & 0.5 & 0.0 & 2.5 & 5.0 & 0.0 \\
    \acrshort{meo} & 0.5 & 1.0 & 1.0 & 0.5 & 0.5 & 0.5 \\
    \bottomrule
    \end{tabular}
    \caption{Initial uncertainties for all test cases}
    \label{tab:initial_sigma}
\end{table}

To better match the bell-shaped Gaussian \acrshort{pdf} and minimize the number of new components generated by the \acrshort{loads} algorithm, the splitting library in \cref{tab:split_lib_coef} with $L=3$ is selected. The \acrlong{ci} is also set to $c=3$. The maximum number of splits and minimum components weight are chosen as $N_{max}=20$ and $\alpha_{min}=10^{-8}$ respectively so that the imposed nonlinearity threshold $\varepsilon_\nu$ is always satisfied. The last is tuned depending on the specific orbit regime and coordinates set as described below.

Keplerian parameters $\bar{\vb*{x}}_{0,kep}$ in \cref{tab:nominal_ics} are converted to cartesian elements $\bar{\vb*{x}}_{0,cart}$ and the initial covariance matrix is obtained as $\vb*{P}_{0,cart} = diag(\vb*{\sigma}^2_{0,cart})$ with $\vb*{\sigma}_{0,cart}$ given by \cref{tab:initial_sigma}. The corresponding \acrshort{da} state $[\vb*{x}_{0,cart}]$ is then initialized according to \cref{eqn:da_input}. For simulations carried out in \acrshortpl{oe} other that cartesian parameters, $[\vb*{x}_{0,cart}]$ must be converted to an \acrshort{ic} in the corresponding elements set. This is achieved with \cref{alg:loads_gmm_algorithm} choosing $\vb*{f}$ as the coordinate transformation from cartesian to target \acrshortpl{oe} and a nonlinearity threshold $\varepsilon_\nu=0.01$. For all proposed scenarios, a single Gaussian kernel is sufficient for this map. Each resulting \acrshort{ic} is then propagated for two Keplerian periods of the corresponding nominal orbit.

\subsection{\glsentrylong{up} in Cartesian Parameters}

A first qualitative analysis on the orbit \acrshort{up} algorithm in cartesian parameters is conducted in this section. Iso-probability lines of the propagated state \acrshort{pdf} and reference \acrshort{mc} samples are plotted in \cref{sec:case_heo,sec:case_leo,sec:case_meo} to visually assess the accuracy of the proposed \acrshort{up} method. A quantitative measure of the accuracy of the proposed \acrlong{mf} method is provided in \cref{sec:rmse}.

\subsubsection{\glsentrylong{heo}}\label{sec:case_heo}

\begin{figure}[ht]
	\begin{subfigure}[b]{0.5\textwidth}
		\centering
		\includegraphics[width=0.95\textwidth]{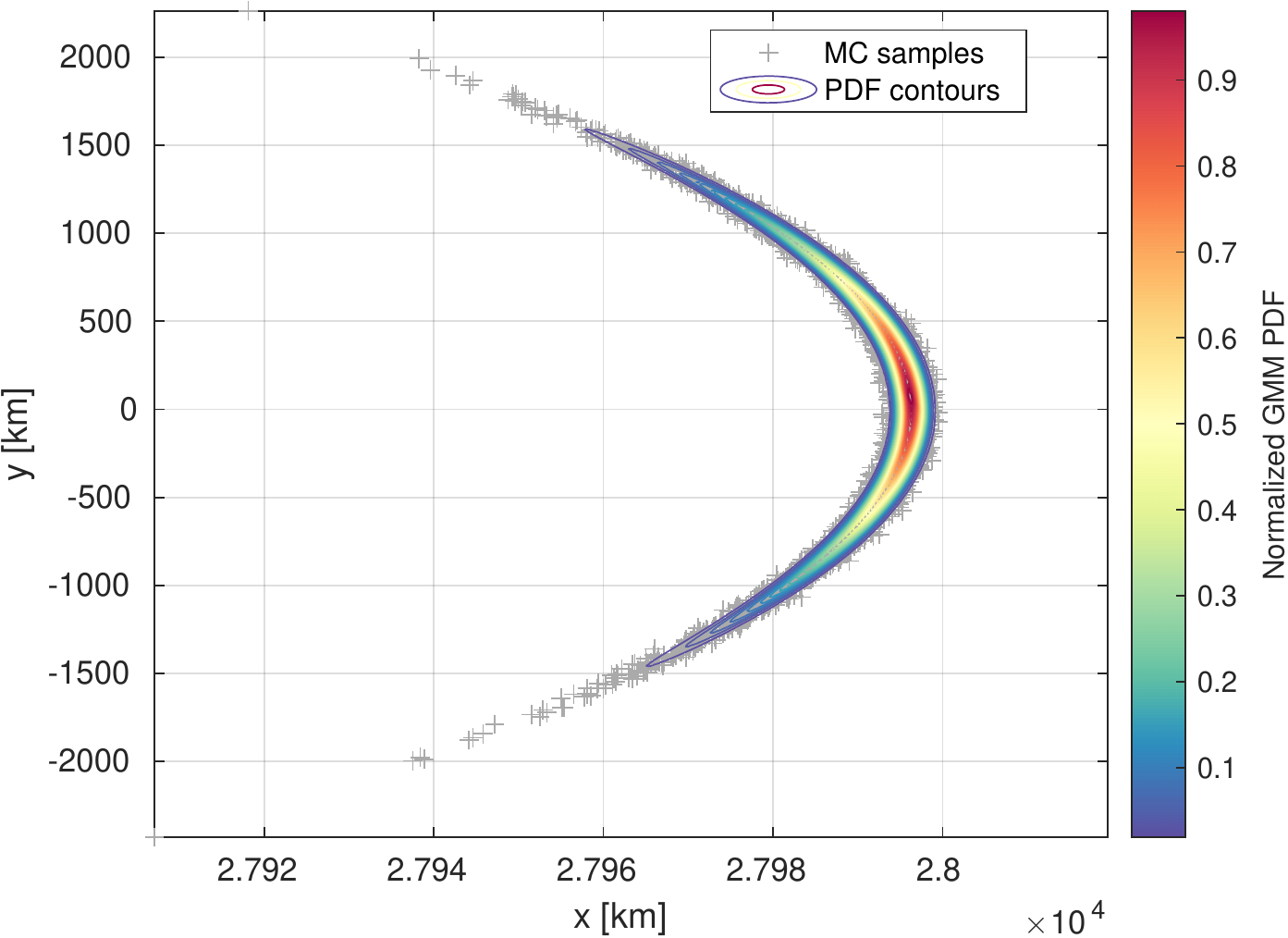}
		\caption{$x-y$ projection}
		\label{fig:heo_xy_contour}
	\end{subfigure}
	\begin{subfigure}[b]{0.5\textwidth}
		\centering
		\includegraphics[width=0.95\textwidth]{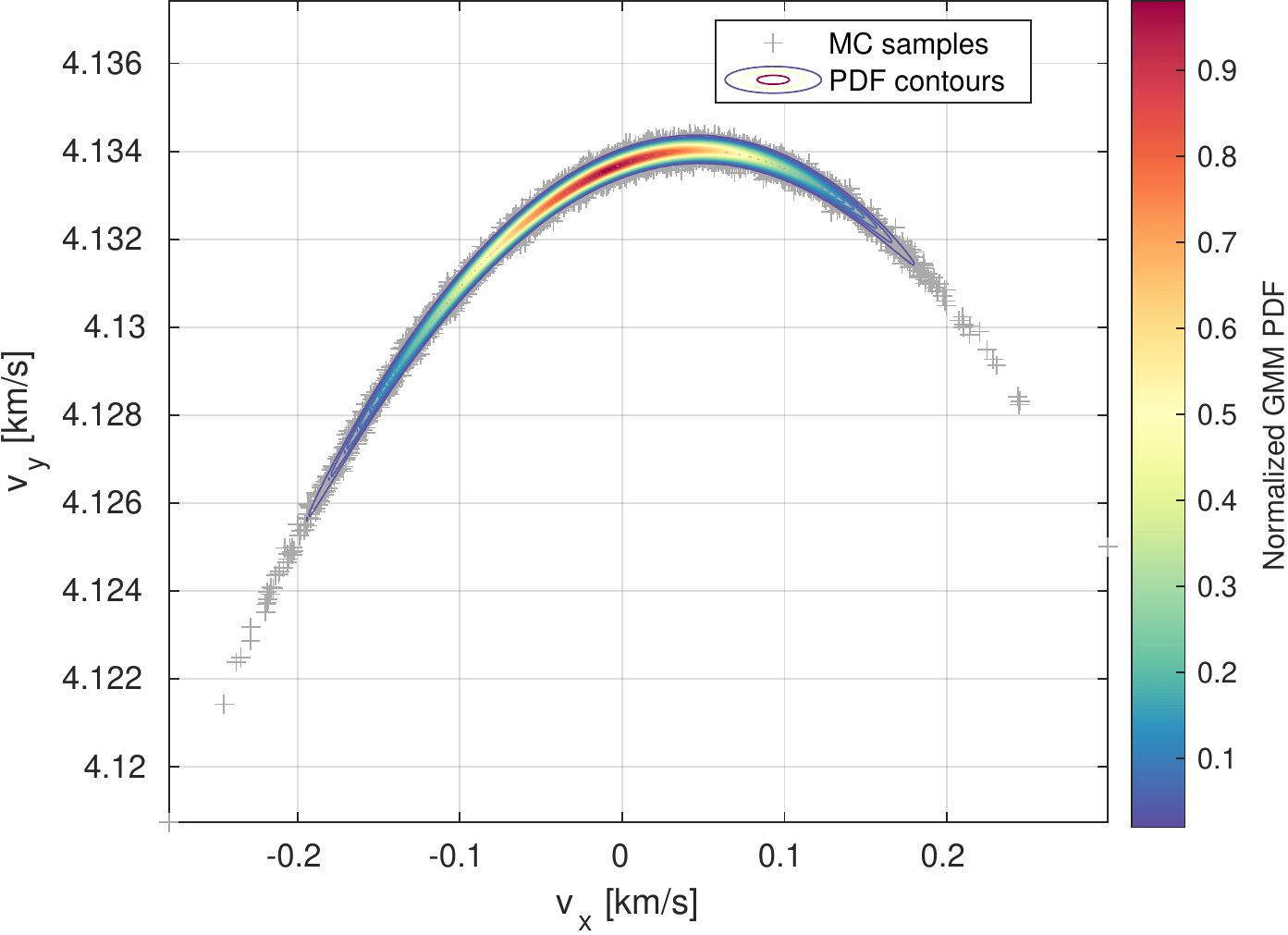}
		\caption{$v_x-v_y$ projection}
		\label{fig:heo_vxvy_contour}
	\end{subfigure}
	\caption{Propagated \glsentryshort{mc} samples and mixture \glsentryshort{pdf} in cartesian parameters for \glsentrylong{heo} case}
	\label{fig:heo_contours}
\end{figure}

Output of the \acrshort{heo} case obtained with $\varepsilon_\nu=0.01$ are displayed in \cref{fig:heo_contours}. \Cref{fig:heo_xy_contour} shows a projection of the propagated uncertainty in position onto the $x-y$ plane while \cref{fig:heo_vxvy_contour} represents the uncertainty in velocity on the same plane. Values of the estimated \acrshort{pdf} are normalized with respect to their maximum located in correspondence of the nominal solution. The final \acrshort{gmm} includes 6561 components to accurately represent the strongly non-Gaussian \acrshort{pdf} as demonstrated by the contour lines following the curvature of the reference \acrshort{mc} samples.

\subsubsection{\glsentrylong{leo}}\label{sec:case_leo}

\begin{figure}[ht]
	\begin{subfigure}[b]{0.5\textwidth}
		\centering
		\includegraphics[width=0.95\textwidth]{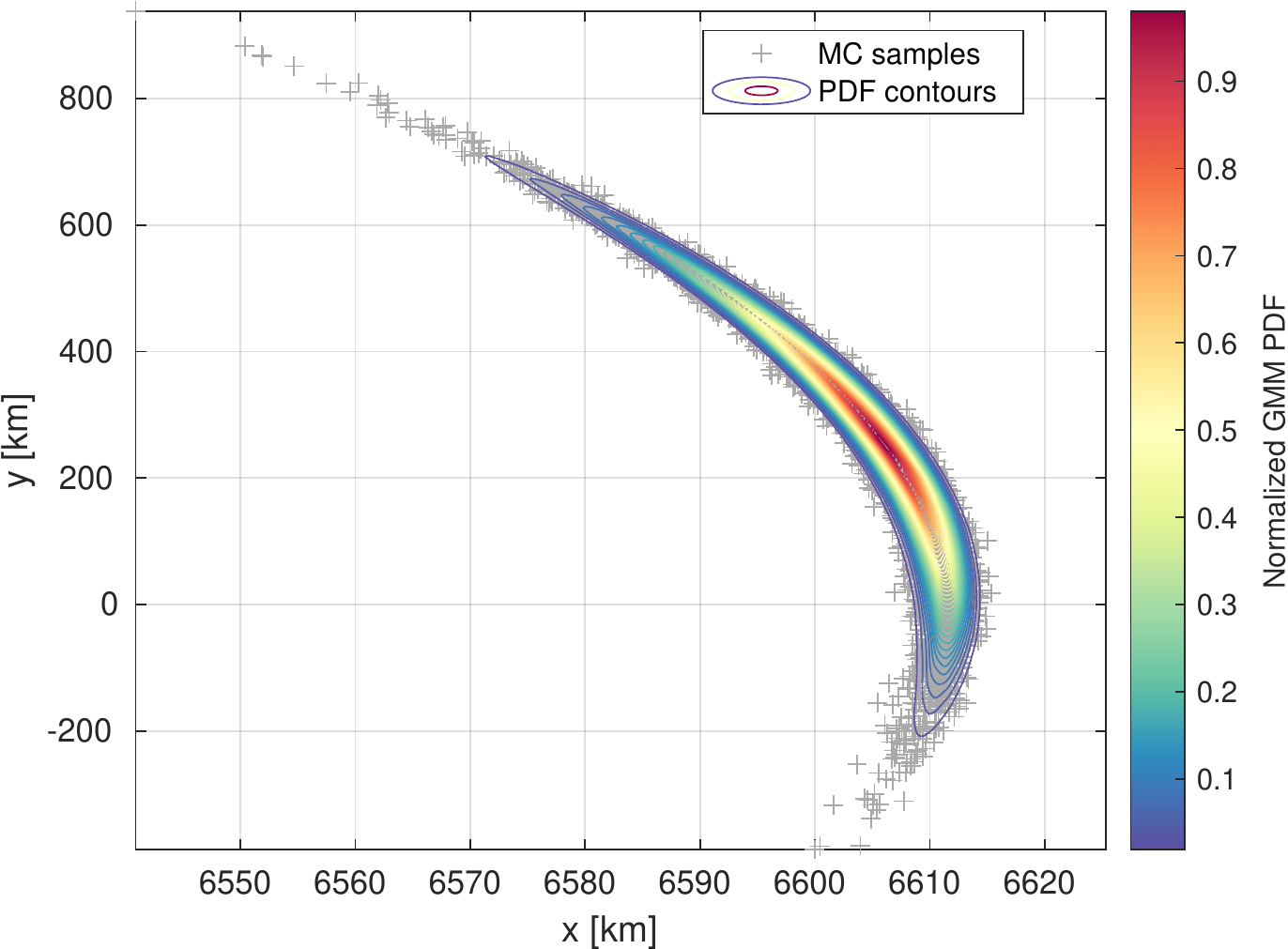}
		\caption{$x-y$ projection}
		\label{fig:leo_xy_contour}
	\end{subfigure}
	\begin{subfigure}[b]{0.5\textwidth}
		\centering
		\includegraphics[width=0.95\textwidth]{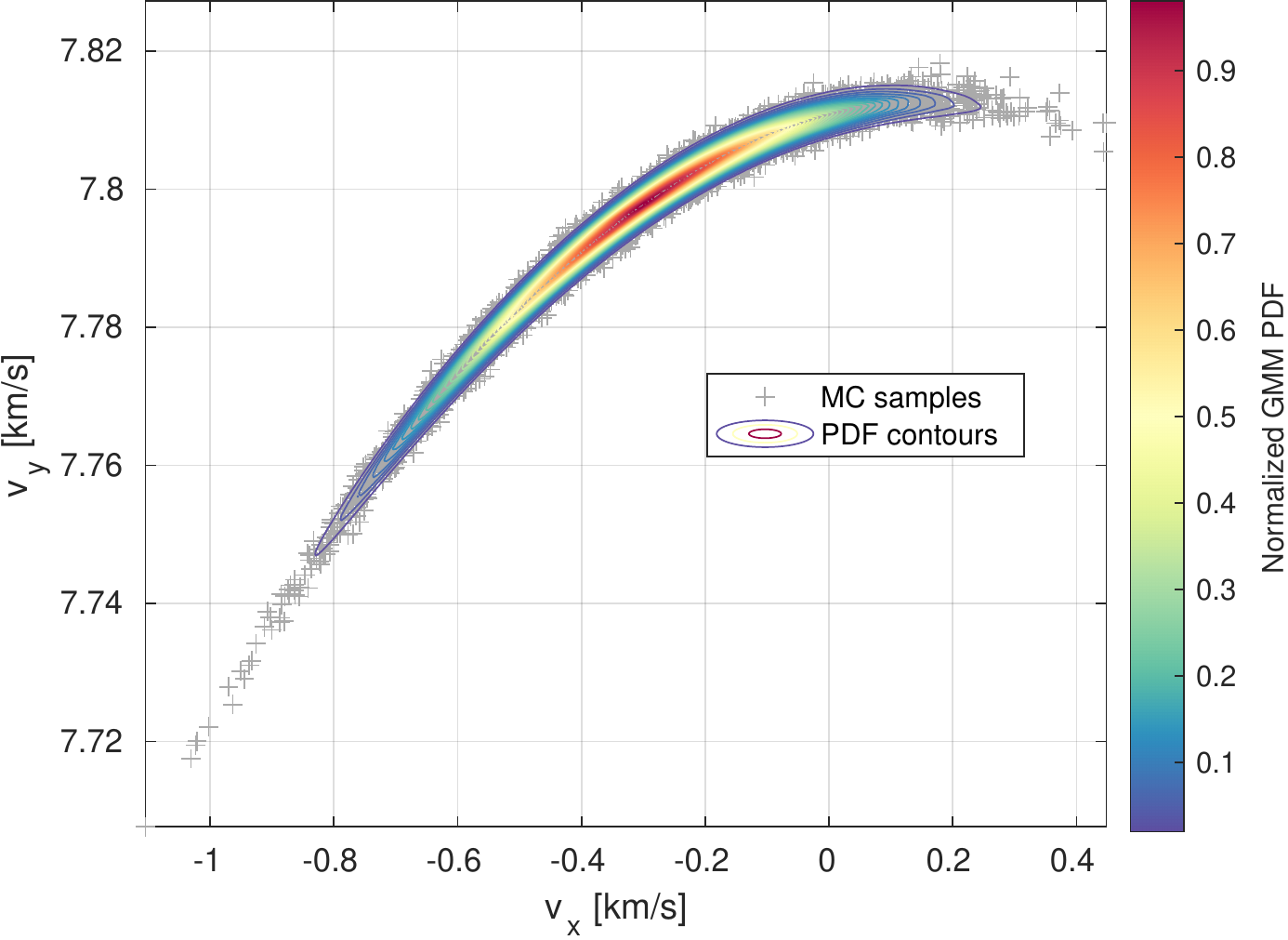}
		\caption{$v_x-v_y$ projection}
		\label{fig:leo_vxvy_contour}
	\end{subfigure}
	\caption{Propagated \glsentryshort{mc} samples and mixture \glsentryshort{pdf} in cartesian parameters for \glsentrylong{leo} case}
	\label{fig:leo_contours}
\end{figure}

The second scenario is the equatorial \acrshort{leo} orbit for which the estimated \acrshort{pdf} is displayed in \cref{fig:leo_contours}. \Cref{fig:leo_xy_contour} shows a projection of the uncertainty in position onto the $x-y$ plane while \cref{fig:leo_vxvy_contour} a projection of the corresponding uncertainty in velocity onto the same plane. As before, the effects of nonlinearity are visible in the curvature of the initially Gaussian \acrshort{mc} samples set, closely followed by the contour lines of the estimated \acrshort{pdf}. Despite the higher nonlinearity threshold $\varepsilon_\nu=0.025$ and the lower eccentricity with respect to the previous case, 2187 kernels are needed to approximate the transformed \acrshort{pdf} mainly due to the lower orbit altitude and the consequent atmospheric drag.

\subsubsection{\glsentrylong{meo}}\label{sec:case_meo}

\begin{figure}[h]
	\begin{subfigure}[b]{0.5\textwidth}
		\centering
		\includegraphics[width=0.95\textwidth]{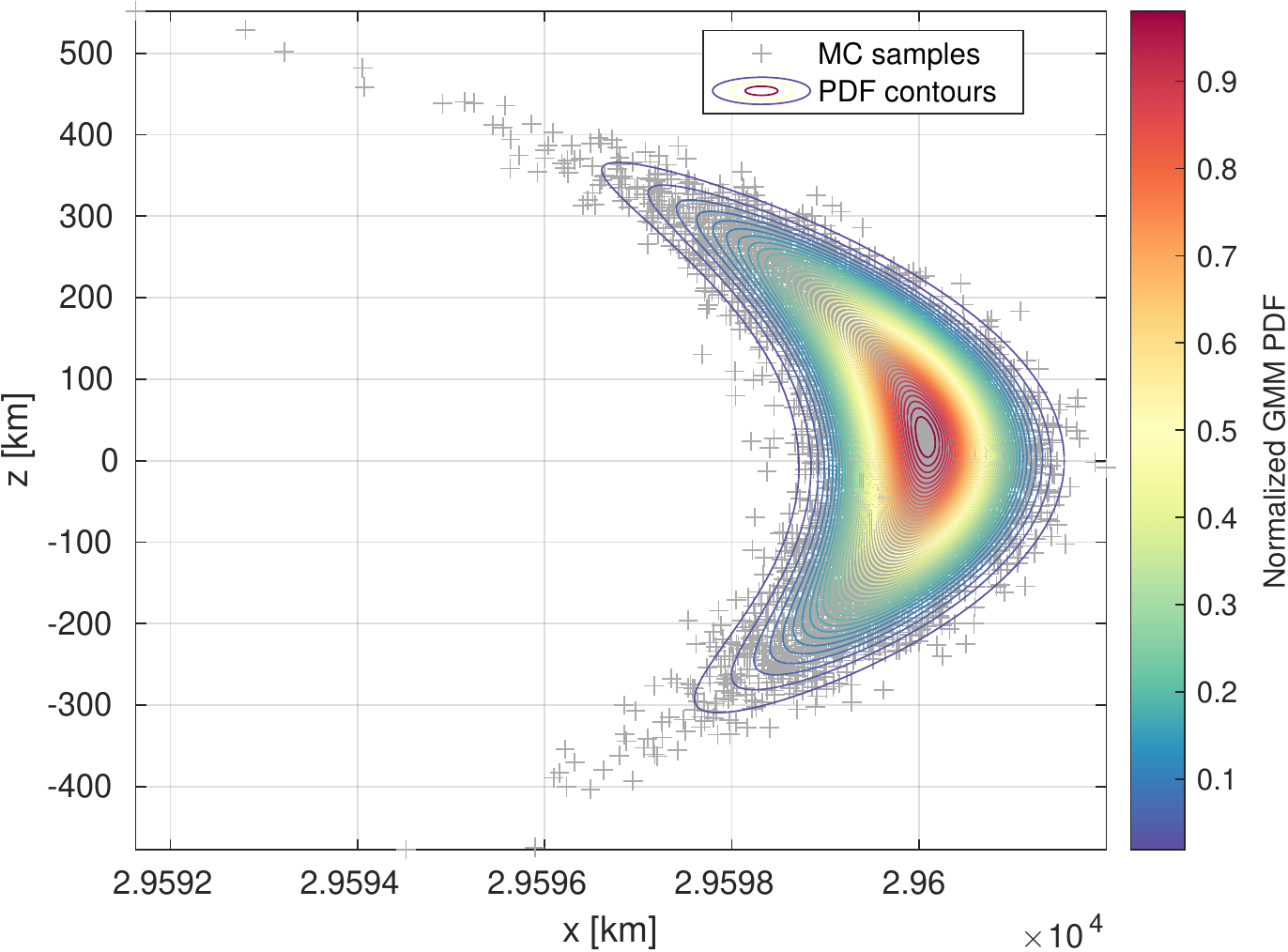}
		\caption{$x-z$ projection}
		\label{fig:meo_xz_contour}
	\end{subfigure}
	\begin{subfigure}[b]{0.5\textwidth}
		\centering
		\includegraphics[width=0.95\textwidth]{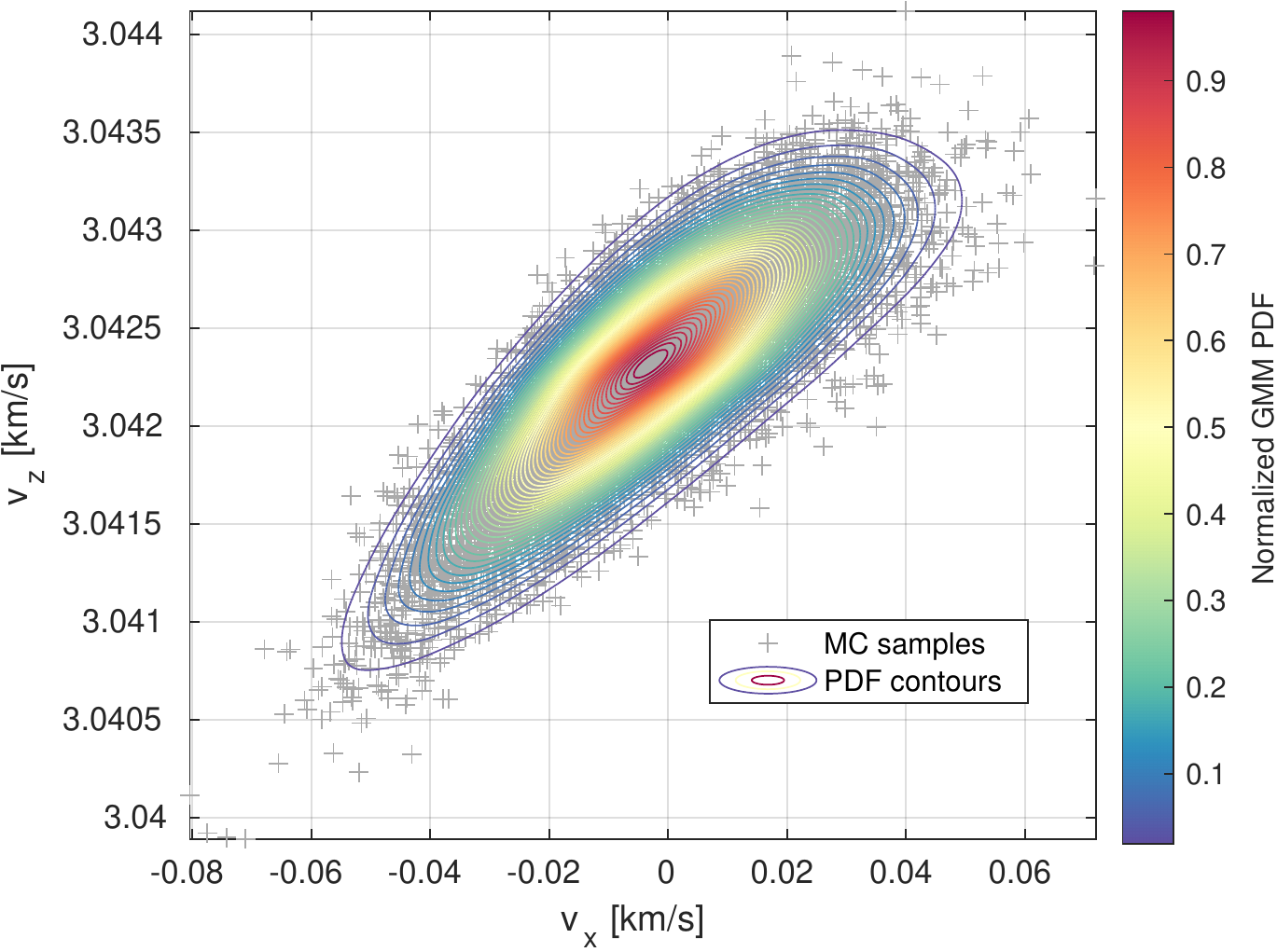}
		\caption{$v_x-v_z$ projection}
		\label{fig:meo_vxvz_contour}
	\end{subfigure}
	\caption{Propagated \glsentryshort{mc} samples and mixture \glsentryshort{pdf} in cartesian parameters for \glsentrylong{meo} case}
	\label{fig:meo_contours}
\end{figure}

The last test case is the circular, inclined orbit of the Galileo constellation. The propagated uncertainties are shown in \cref{fig:meo_contours} together with the reference \acrshort{mc} samples. \Cref{fig:meo_xz_contour} represents the uncertainty in position projected on the $x-z$ plane while \cref{fig:meo_vxvz_contour} the corresponding uncertainty in velocity. In this case, setting $\varepsilon_\nu=0.01$ leads to 729 \acrshort{gmm} components required to approximate the propagated \acrshort{pdf}. The great reduction in the number of kernels with respect to the previous cases is mainly due to the lower eccentricity (compared to the \acrshort{heo} case), the longer orbital period and the absence of atmospheric drag (compared to the \acrshort{leo} case). The lower number of components is thus a hint of the weaker nonlinearity of the dynamics in this scenario as demonstrated by the shape of the contour lines in \cref{fig:meo_vxvz_contour}. This behavior also highlights the advantages of an adaptive splitting scheme with respect to more conventional \acrshortpl{gmm} with fixed number of components: the final number of Gaussian kernels is guaranteed to be the minimum required to satisfy the imposed nonlinearity threshold.

\subsection{Sensitivity to State Parametrization and \glsentryshort{loads}-\glsentryshort{gmm} Parameters}\label{sec:sensitivity}

\subsubsection{Sensitivity to State Parametrization}\label{sec:sensitivity_oe}

\begin{figure}[ht]
    \centering
    \includegraphics[width=0.5\textwidth]{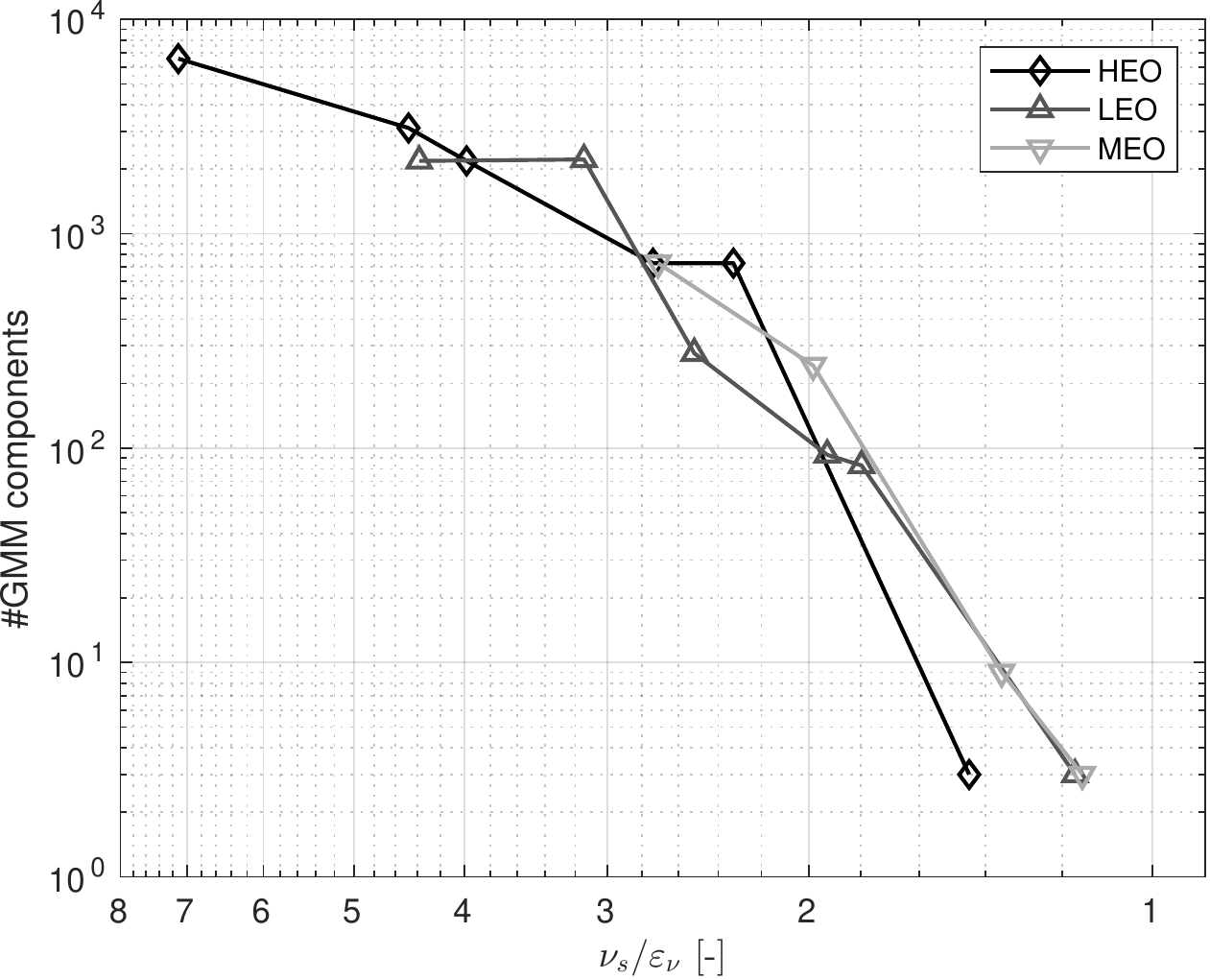}
    \caption{Number of \acrshort{gmm} components as function of normalized single-domain \glsentryshort{nli} $\nu_{s}/\varepsilon_\nu$}
    \label{fig:nb_kernels_vs_nu_single}
\end{figure}

This section discusses the impact of the state parametrization, i.e. the selected coordinates set, on the overall \acrshort{mf} algorithm. Twenty-one simulations were carried out. They correspond to all possible combinations among orbit regime (i.e. \acrshort{heo}/\acrshort{leo}/\acrshort{meo}), coordinates set (i.e. cartesian/equinoctial/\acrshort{mee}/alternate elements) and position angles (i.e. $L,\lambda$, not applicable for cartesian parameters). For each case, the imposed nonlinearity threshold $\varepsilon_\nu$, final number of Gaussian kernels $N_f$ and single-domain \acrshort{nli} $\nu_s$ are reported in \cref{tab:nu_single_nb_kernels}. $\nu_s$ corresponds to the \acrshort{nli} of the final domain $[\vb*{x}(t)]$ when mapped without \acrshort{ads} and a direct measure of nonlinearity of the different representations. For a fixed threshold $\varepsilon_\nu$, $\nu_s$ is thus expected to be directly correlated with the final number of components $N_f$. This fact is demonstrated in \cref{fig:nb_kernels_vs_nu_single} in which $N_f$ is plotted as function of the normalized single-domain \acrshort{nli} $\nu_s/\varepsilon_\nu$ (only the cases for which $\nu_s > \varepsilon_\nu$ and thus $N_f>1$ are shown in the picture). For each orbit regime, the number of components decreases monotonically for decreasing $\nu_s/\varepsilon_\nu$ as expected by construction of the \acrshort{loads} algorithm.

\begin{table}[ht]
\centering
\small
\sisetup{round-mode=places,round-precision=3}
\begin{tabular}{rc S[scientific-notation=true,table-format=1.3e-2] cc S[scientific-notation=true,table-format=1.3e-2] cc S[scientific-notation=true,table-format=1.3e-2] cc}
\toprule
& & \multicolumn{3}{c}{\acrlong{heo}} & \multicolumn{3}{c}{\acrlong{leo}} & \multicolumn{3}{c}{\acrlong{meo}} \\ \midrule
\acrshortpl{oe} & & $\nu_s$ & $\varepsilon_\nu$ & $N_f$ & $\nu_s$ & $\varepsilon_\nu$ & $N_f$ & $\nu_s$ & $\varepsilon_\nu$ & $N_f$ \\ \midrule
Cartesian & $-$ & 0.0712063792682708 & 0.01  & 6561 & 0.1095655792679876 & 0.025 & 2187 & 0.0270885017785335 & 0.01  & 729 \\
Equinoctial & $L$ & 0.0134340767596835 & 0.003 & 3115 & 0.0251826657009531 & 0.01  & 277  & 0.0040672823614179 & 0.003 & 9   \\
\acrshortpl{mee} & $L$ & 0.0119511344382026 & 0.003 & 2187 & 0.0314460477883771 & 0.01  & 2221 & 0.0059464064967170 & 0.003 & 243 \\
Alternate & $L$ & 0.0082101011805195 & 0.003 & 729  & 0.0192651487099065 & 0.01  & 93   & 0.0027969071275029 & 0.003 & 1   \\
Equinoctial & $\lambda$ & 0.0043420514522953 & 0.003 & 3    & 0.0116954699553937 & 0.01  & 3    & 0.0015764072902437 & 0.003 & 1   \\
\acrshortpl{mee} & $\lambda$ & 0.0069824346210532 & 0.003 & 729  & 0.0179780982246813 & 0.01  & 83   & 0.0034572238627005 & 0.003 & 3   \\
Alternate & $\lambda$ & 0.0000071512744626 & 0.003 & 1    & 0.0037825961778005 & 0.01  & 1    & 0.0000173655491134 & 0.003 & 1   \\
\bottomrule
\end{tabular}
\caption{Single-domain \glsentryshort{nli} $\nu_{s}$, nonlinearity threshold $\varepsilon_\nu$ and number of \glsentryshort{gmm} components for all tests cases}
\label{tab:nu_single_nb_kernels}
\end{table}

From \cref{tab:nu_single_nb_kernels} it is clear that cartesian parameters are the most nonlinear since they lead to larger number of kernels despite the higher thresholds selected in this case. Moreover, the use of alternate elements and mean longitude $\lambda$ results in a quasi-linear transformation that is correctly captured by a single polynomial as seen in the last row of \cref{tab:nu_single_nb_kernels}.

Unfortunately, coordinates sets can rarely be chosen freely. Many applications such as collision probability estimation require the state uncertainty at the \acrfull{tca} expressed in cartesian parameters. Adding a coordinate transformation at the end of the \acrshort{loads} target function $\vb*{f}$ will be detrimental in this case since the \acrshort{nli} will still be computed on cartesian parameters, thus loosing all benefits of alternative state parametrizations. However, this is only true for analytical propagators such as \acrshort{sgp}4 which require a single function evaluation to compute the final state. Conversely, for numerical propagators splits are allowed at any time along the trajectory. Since each split will then correspond to $L-1$ new trajectories to be propagated\footnote{\label{ftn:loads_num_prop}details on how numerical integrators can be embedded in the \acrshort{loads} framework are provided in \cite{Losacco}}, it is of great interest to delay as much as possible these splits and propagate a single expansion for the entire time span. This behavior could be achieved using alternate equinoctial elements for the propagation itself followed by a coordinate transformation of the single final state in which all splits will occur.

\subsubsection{Sensitivity to Minimum Kernel Weight}\label{sec:alpha_sensitivity}

Since the number of \acrshort{gmm} components output of the \acrshort{mf} method are bounded by the maximum allowed number of splits $N_{max}$ and the minimum kernel weight $\alpha_{min}$, it is of interest to analyze the impact of $\alpha_{min}$ on the accuracy of the estimated \acrshort{pdf}. This analysis is conducted in cartesian parameters since this coordinates set generates the highest number of splits and is thus the one which is penalized the most by large thresholds on $\alpha_{min}$. For the three test cases presented above, multiple solutions were obtained for different values of $\alpha_{min}\in[1, 10^{-8}]$ and the corresponding \acrshort{lam} with respect to the same set of \acrshort{mc} samples computed using \cref{eqn:lam_dmm}. These values were then normalized with respect to the solution with the highest allowed number of kernels, i.e.

\begin{equation}
    \mathcal{L}_{n,MF}(\alpha_{min}) = \frac{\mathcal{L}_{MF}(\alpha_{min})}{\mathcal{L}_{MF}(\alpha_{min})\rvert_{\alpha_{min}=10^{-8}}}
    \label{eqn:normalized_lam_mf}
\end{equation}

\begin{figure}[ht]
    \begin{subfigure}[b]{0.5\textwidth}
        \centering
        \includegraphics[width=0.95\textwidth]{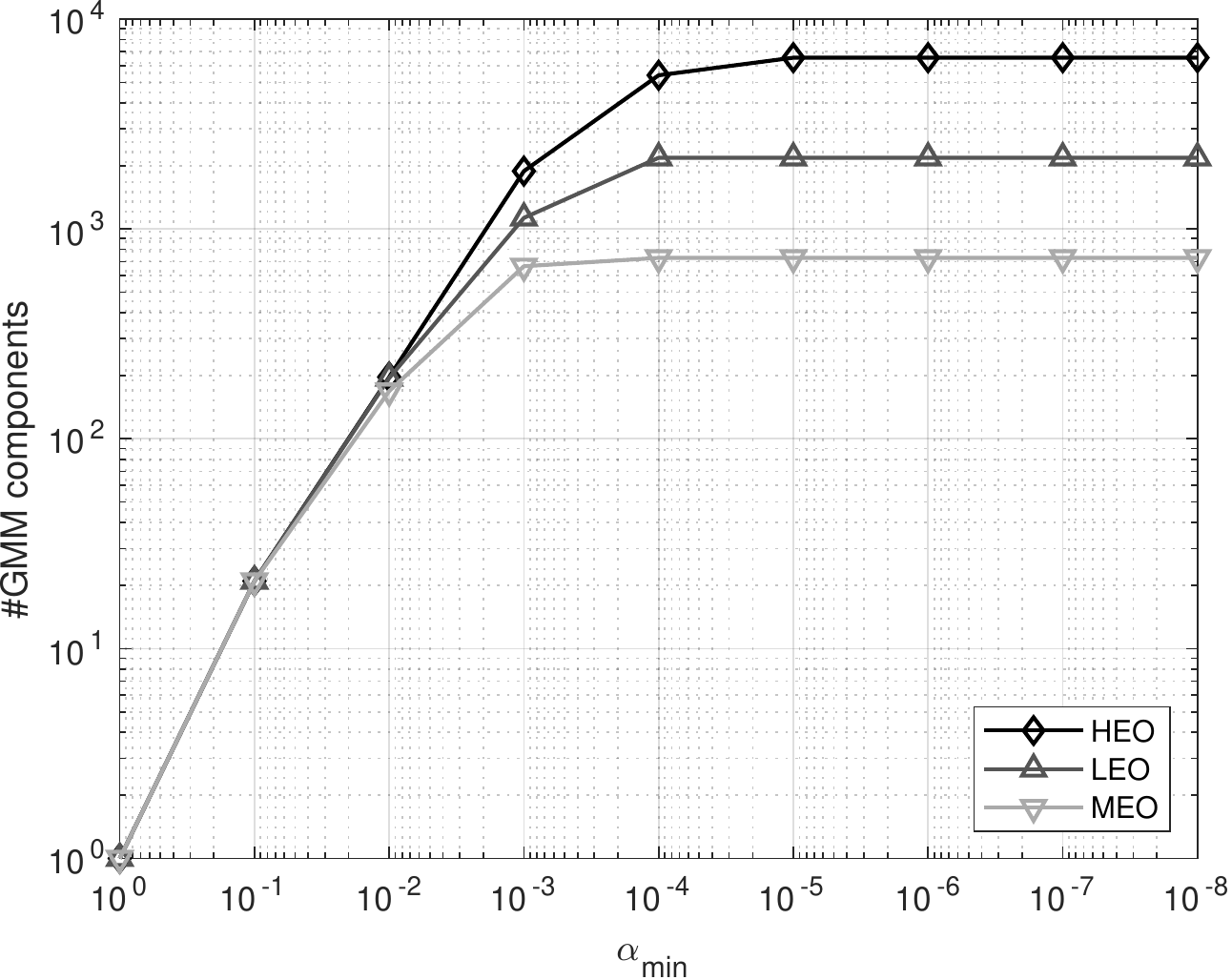}
        \caption{Number of \glsentryshort{gmm} components}
        \label{fig:nb_comps_vs_alpha}
    \end{subfigure}
    \begin{subfigure}[b]{0.5\textwidth}
        \centering
        \includegraphics[width=0.95\textwidth]{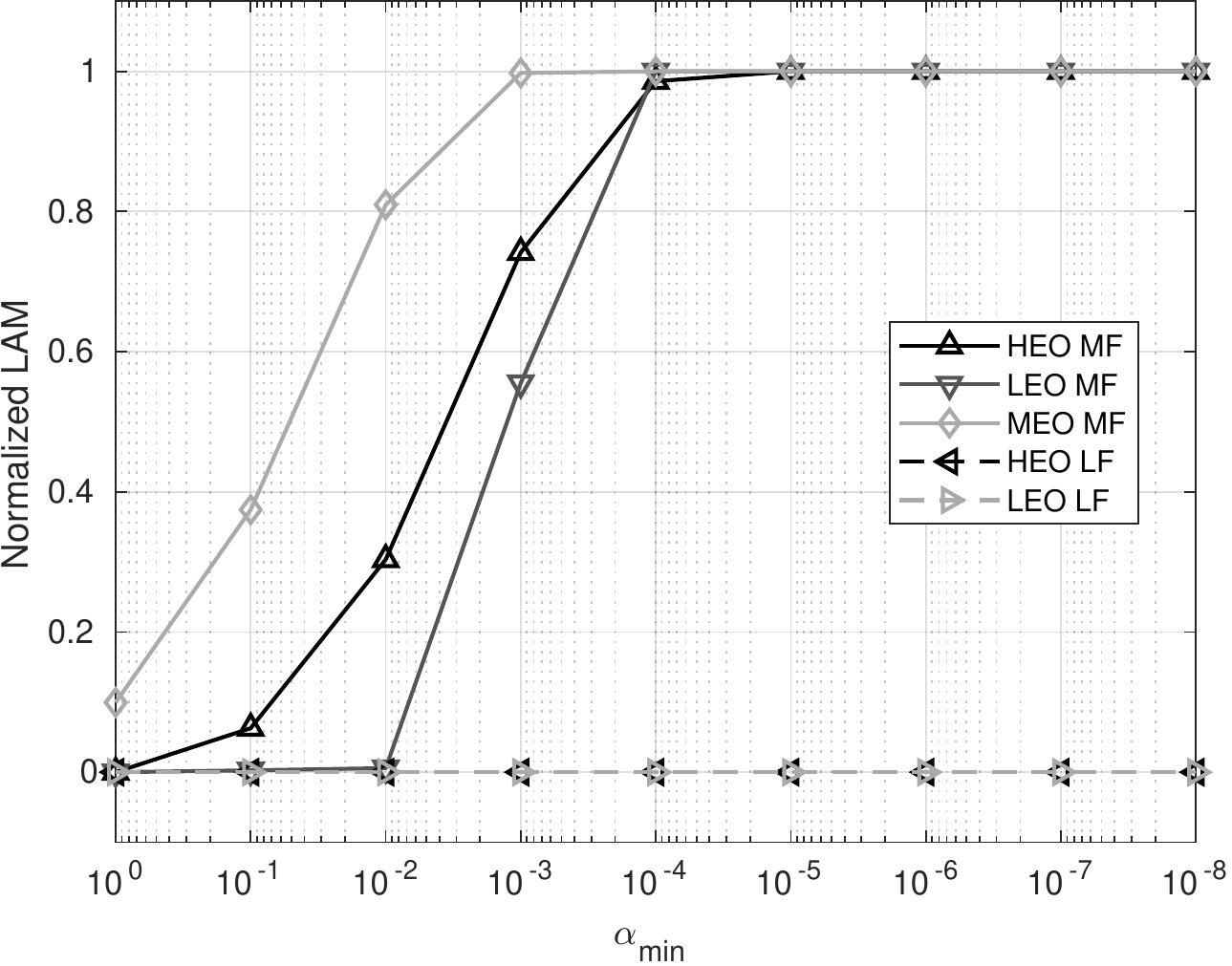}
        \caption{Normalized \glsentryshort{lam}}
        \label{fig:lam_vs_alpha}
    \end{subfigure}
    \caption{Number of \glsentryshort{gmm} components and normalized \glsentryshort{lam} as function of $\alpha_{min}$}
    \label{fig:comps_lam_vs_alpha}
\end{figure}

The final number of components $N_f$ and normalized \acrshort{lam} $\mathcal{L}_n$ are shown in \cref{fig:comps_lam_vs_alpha} as a function of $\alpha_{min}$.

As can be seen in \cref{fig:nb_comps_vs_alpha}, for the \acrlong{heo} test case no more than 6561 kernels with $\alpha_i>10^{-5}\ \forall i$ are needed to satisfy the imposed nonlinearity threshold $\varepsilon_\nu$. For $\alpha_{min}<10^{-5}$ a plateau is in fact reached in both figures, thus signifying the same \acrshort{pdf} estimation is obtained for smaller values. On the contrary, looking at \cref{fig:lam_vs_alpha} from right to left a rapid decrease in $\mathcal{L}_n$ is observed for $\alpha_{min}\to1$, i.e. for a number of \acrshort{gmm} components that is progressively reduced to a single Gaussian kernel. Similar conclusions can be drawn for the \acrshort{leo} and \acrshort{meo} cases in which the plateau is reached for $\alpha_{min}=10^{-4}$ corresponding to 2187 and 729 components respectively. These results demonstrate the flexibility and accuracy of the proposed \acrshort{mf} approach which eliminates the needs of a priori estimation of the required number of kernels to guarantee a certain accuracy thus minimizing the computational effort given the nonlinearity content of the problem at hand.\par

\subsection{Accuracy of the \glsentrytitlecase{mf}{long} Method}\label{sec:rmse}

The \acrfull{rmse} is introduced in this section to quantify the error made by the \acrlong{mf} method with respect to a reference \acrshort{mc} simulation. The \acrshort{rmse} is defined as

\begin{equation}
    \vb*{e}_{RMSE} = \sqrt{\dfrac{1}{N} \sum_{i=1}^N \left(\vb*{x}_i-\hat{\vb*{x}}_i\right)^2}
    \label{eqn:rmse_def}
\end{equation}

with $\{\vb*{x}_i\},\{\hat{\vb*{x}}_i\}$ expected and actual samples respectively. Note that all operations in \cref{eqn:rmse_def} are carried out element-wise so that $\vb*{e}_{RMSE}$ has the same dimension of $\hat{\vb*{x}}_i,\vb*{x}_i\ \forall i$.

To compute $\{\vb*{x}_n\},\{\hat{\vb*{x}}_i\}$ in the case of orbit \acrshort{up}, $N$ random samples $\{\vb*{x}_i(t_0)\}$ are firstly drawn from the initial distribution $p_{\vb*{X}}(\vb*{x}(t_0))$ in \cref{eqn:initial_normal_distribution_mf}. $\{\vb*{x}_i\}$ is then simply obtained with a \acrlong{pw} propagation of $\{\vb*{x}_i(t_0)\}$ in \acrshort{hf} dynamics. Conversely, $\{\hat{\vb*{x}}_i\}$ is computed evaluating the Taylor polynomials $[\vb*{x}_{MF}^{(p)}(t)]$ in $\{\vb*{x}_i(t_0)\}$ thus avoiding the needs for additional propagation steps. For each entry in \cref{tab:nu_single_nb_kernels}, three samples sets are computed. They correspond to the \acrshort{lf} solution and the two \acrlong{mf} ones obtained with a polynomial shift applied in either osculating or \acrshort{tle} elements space using \cref{eqn:mf_expansion_splitted} or \cref{alg:mf_shift_tle} respectively. For the shake of clarity, only a subset of the computed errors is reported in \cref{tab:rmse_summary}.

\begin{table}[ht]
    \centering
    \small
    \sisetup{round-mode=places,round-precision=3}
    \begin{tabular}{cl *6{S[scientific-notation=true,table-format=1.3e-2]}}
    \toprule
    \multirow{8}{*}{\acrshort{heo}} & & {$x$, \si{\km}} & {$y$, \si{\km}} & {$z$, \si{\km}} & {$v_x$, \si{\km\per\s}} & {$v_y$, \si{\km\per\s}} & {$v_z$, \si{\km\per\s}} \\ \cmidrule{2-8}
    & \acrshort{lf} & 0.836346374734357    & 22.9340274537236     & 0.41545450119153     & 0.00283860221373275  & 0.000102235246382906 & 9.12647380902107e-05 \\
    & \acrshort{mf}\textsubscript{osc} & 0.0242473505041237   & 0.00785463918116089  & 0.000665698295298176 & 8.80377950898014e-07 & 3.61717891522195e-06 & 1.34409484710056e-07 \\
    & \acrshort{mf}\textsubscript{TLE} & 0.000697723079789579 & 0.0093453915285617   & 0.000379199894605152 & 9.77065946522633e-07 & 1.89072419275017e-06 & 1.00869145732918e-07 \\ \cmidrule{2-8}
    & & {$a$, \si{\km}} & {$f$} & {$g$} & {$h$} & {$k$} & {$L$, \si{\radian}} \\ \cmidrule{2-8}
    & \acrshort{lf} & 0.0576557886832245   & 2.10884230308582e-05 & 5.02829793290442e-06 & 1.10750854378778e-05 & 7.37784365959453e-06 & 0.000819301738233098 \\
    & \acrshort{mf}\textsubscript{osc} & 0.000538602426110117 & 1.92614209458474e-08 & 3.0865072991766e-08  & 1.17960210913868e-08 & 6.61549651600255e-09 & 1.21609675925911e-06 \\
    & \acrshort{mf}\textsubscript{TLE} & 0.048043351202763    & 1.09833981590107e-06 & 5.84163472953534e-08 & 1.22796367333396e-08 & 6.56952208691966e-09 & 1.04189201779176e-06 \\ \midrule
    \multirow{8}{*}{\acrshort{leo}} & & {$p$, \si{\km}} & {$f$} & {$g$} & {$h$} & {$k$} & {$L$, \si{\radian}} \\ \cmidrule{2-8}
    & \acrshort{lf} & 0.0307887638175996   & 2.95632247467367e-05 & 2.81397044865738e-05 & 1.2407504571065e-05  & 4.62151585722597e-06 & 0.000417535892667894 \\
    & \acrshort{mf}\textsubscript{osc} & 0.000509614726516381 & 5.82847899469287e-08 & 1.40879557130554e-07 & 8.1005451223284e-09  & 1.66144849003934e-09 & 7.07281507953398e-06 \\
    & \acrshort{mf}\textsubscript{TLE} & 0.0060777327046681   & 1.44190061935856e-06 & 1.1434863677771e-06  & 8.1233543695718e-09  & 1.64248591119788e-09 & 8.42911649310957e-06 \\ \cmidrule{2-8}
    & & {$n$, \si{\radian\per\s}} & {$f$} & {$g$} & {$h$} & {$k$} & {$\lambda$, \si{\radian}} \\ \cmidrule{2-8}
    & \acrshort{lf} & 8.95614138853844e-09 & 2.95346721138355e-05 & 2.81334613645766e-05 & 1.24060400499539e-05 & 4.62115448764967e-06 & 0.00046239352815581  \\
    & \acrshort{mf}\textsubscript{osc} & 6.89756660845075e-10 & 2.94676935770315e-07 & 1.06839376414294e-06 & 1.6245438387932e-07  & 2.90295290118332e-08 & 7.52430893855102e-05 \\
    & \acrshort{mf}\textsubscript{TLE} & 1.52431258697598e-09 & 1.41044392299758e-06 & 4.13161239628592e-06 & 1.63017794255407e-07 & 2.84484752755287e-08 & 7.52512349145657e-05 \\ \midrule
    \multirow{8}{*}{\acrshort{meo}} & & {$x$, \si{\km}} & {$y$, \si{\km}} & {$z$, \si{\km}} & {$v_x$, \si{\km\per\s}} & {$v_y$, \si{\km\per\s}} & {$v_z$, \si{\km\per\s}} \\ \cmidrule{2-8}
    & \acrshort{lf} & 0.031571015773124    & 1.89849326619983     & 2.80021725828257     & 0.000419724278574856 & 2.79680008528236e-06 & 1.15401711670443e-05 \\
    & \acrshort{mf}\textsubscript{osc} & 0.00274342936393501  & 0.00055121986959116  & 0.000813690171944583 & 1.30695728678572e-07 & 1.65756222088673e-07 & 2.53386436752794e-07 \\
    & \acrshort{mf}\textsubscript{TLE} & 0.0508923105272844   & 0.0318254127358066   & 0.0472510382431536   & 3.56140841759247e-06 & 3.69198581702709e-06 & 5.46996631647552e-06 \\ \cmidrule{2-8}
    & & {$n$, \si{\radian\per\s}} & {$f$} & {$g$} & {$h$} & {$k$} & {$L$, \si{\radian}} \\ \cmidrule{2-8}
    & \acrshort{lf} & 7.11536463236826e-10 & 4.79960745259174e-06 & 8.64513664966206e-08 & 7.06607175253785e-07 & 1.76948087432167e-07 & 0.000114411806049457 \\
    & \acrshort{mf}\textsubscript{osc} & 1.77226502986563e-11 & 9.85958830896377e-08 & 1.43853458746044e-08 & 7.86367355016925e-09 & 1.05360379473876e-08 & 3.97568160859861e-07 \\
    & \acrshort{mf}\textsubscript{TLE} & 1.70630691338012e-10 & 3.33687748143696e-06 & 6.40709604889585e-06 & 7.83594589452714e-09 & 1.04933581002135e-08 & 1.32554866902391e-05 \\
    \bottomrule
    \end{tabular}
    \caption{\glsentryshort{rmse} error for various test cases}
    \label{tab:rmse_summary}
\end{table}

Data in \cref{tab:rmse_summary} demonstrate the accuracy of the proposed \acrshort{up} algorithm in terms of both absolute \acrshort{rmse} with respect to the reference \acrshort{mc} simulation and relative error with respect to the \acrshort{lf} counterpart. When using cartesian parameters, the absolute \acrshort{rmse} in position is always below $\SI{10}{\m}$, while the errors in semimajor axis $a$ and semilatus rectum $p$ are smaller than $\SI{1}{\m}$ if the true longitude $L$ is used and below $\SI{10}{\m}$ otherwise. \acrshortpl{rmse} in $f,g,h,k$ are below $\num{1e-6}$ for all types of \acrshortpl{oe} with \acrshortpl{mee} performing best if combined with $L$. Finally, errors in mean and true longitudes $\lambda,L$ exhibit a stronger dependence on the orbit regime, with \acrshortpl{rmse} that vary from $\SI{1e-5}{\radian}$ in \acrshort{leo} to $\SI{1e-7}{\radian}$ in \acrshort{meo}. The \acrlong{mf} correction is key to achieve such small errors since it guarantees a major improvement over the \acrshort{lf} solution. In particular, the \acrshort{rmse} is reduced by up to four \acrfullpl{oom} for propagations in cartesian parameters, with the polynomial shift better performed in osculating elements space. The only exception is the $x$ coordinate in \acrshort{heo} regime, for which a shift in \acrshort{tle} space further reduces the \acrshort{rmse} by two \acrshortpl{oom}. Concerning the semimajor axis $a$, semilatus rectum $p$ and mean motion $n$, a shift in osculating elements space is always preferred while the application of \cref{alg:mf_shift_tle} could result in \acrshortpl{rmse} close to the \acrshort{lf} solution. This behavior is more marked in \acrshort{heo} and \acrshort{meo} regimes. Similar arguments apply to $f$ for which a direct shift results in a two \acrshortpl{oom} improvement with respect to the \acrshort{lf} baseline. Application of \cref{alg:mf_shift_tle} results in only one \acrshort{oom} improvement for the \acrshort{heo} and \acrshort{leo} cases and no improvement in \acrshort{meo}. For the $g$ coordinate, the \acrshort{rmse} of the \acrlong{mf} solution is one to two \acrshortpl{oom} smaller that that of the \acrshort{lf} counterpart if the shift is performed in osculating elements space. Conversely, the accuracy is similar in \acrshort{heo} while it degrades in \acrshort{leo} and especially in \acrshort{meo} for shifts in \acrshort{tle} space. Contrarily to the previous cases, the accuracy in $h,k$ does not depend on the space in which the shift is performed and both \acrlong{mf} solutions guarantee a two to three \acrshortpl{oom} improvement with respect to the \acrshort{lf} one. The same behavior is observed for the mean longitude $\lambda$ and the true longitude $L$ excluding the \acrshort{meo} regime for the last. In this case, the error is one to three \acrshortpl{oom} smaller after the correction is applied. Finally, application of \cref{alg:mf_shift_tle} in \acrshort{meo} regime worsens the \acrshort{rmse} of the \acrshort{lf} solution while a three \acrshortpl{oom} improvement is otherwise obtained. To summarize, the accuracy of the proposed \acrlong{mf} \acrshort{up} method is within few meters in position over a wide range of orbit regimes. Moreover, a simple polynomial shift applied in osculating elements space through \cref{eqn:mf_expansion_splitted} is preferred to \cref{alg:mf_shift_tle} since the last could result in virtually no improvements over the baseline \acrshort{lf} solution.

\subsection{Computational Savings}\label{sec:computational_load}

The \acrshort{mf} method presented in this paper has been introduced as an accurate and computationally efficient alternative to more conventional \acrshort{hf} approaches for orbit \acrshort{up}. In this section, the computational load of the proposed method is compared to that of an equivalent \acrshort{hf} solution. For this purpose, a \acrshort{da} version of the \acrshort{hf} numerical propagator is used with the combined \acrshort{loads}-\acrshort{gmm} algorithm to propagate the initial uncertainty directly in \acrshort{hf} thus avoiding the needs for a posteriori correction\textsuperscript{\ref{ftn:loads_num_prop}}. All numerical simulations were run on an Intel Core i5-9400H @ \SI{2.50}{\giga\hertz} with 16GB of RAM running Ubuntu 18.04 LTS. The total runtime was approximated from that of a single propagation as

\begin{subequations}
    \begin{align}
        T_{MF} &\approx L \cdot \left(t_{DA,LF} + t_{PW,HF} \right) \\
        T_{HF} &\approx L \cdot t_{DA,HF}
    \end{align}
\end{subequations}

where $L$ is the number of Gaussian kernels and the subscripts $DA$ and $PW$ indicate whether the propagation is performed in \acrshort{da} or \acrlong{pw}. The runtimes for single kernels in cartesian parameters are summarized in \cref{tab:runtime}.

\begin{table}[ht]
    \centering
    \begin{tabular}{cccccc}
        \toprule
        Regime & $t_{DA,LF},\si{\milli\s}$ & $t_{PW,HF},\si{\milli\s}$ & $t_{MF},\si{\milli\s}$ & $t_{DA,HF},\si{\milli\s}$ & $t_{DA,HF}/t_{MF},-$ \\ \midrule 
        \acrshort{heo} & 2.805 & 67.395 & 70.200 & 1442.065 & 20.542 \\
        \acrshort{leo} & 1.855 & 353.120 & 354.975 & 5389.915 & 15.184 \\
        \acrshort{meo} & 2.350 & 60.560 & 62.910 & 1176.670 & 18.704 \\
        \bottomrule
    \end{tabular}
    \caption{Runtime for a single propagation in cartesian parameters for all test cases}
    \label{tab:runtime}
\end{table}

As seen in \cref{tab:runtime}, the \acrshort{mf} method is 15 to 20 times faster than the \acrshort{hf} counterpart with the \acrlong{heo} case providing the highest speedup. Moreover, within the \acrshort{mf} procedure most of the time is spent in the correction step which must be carried out numerically for maximum accuracy. This fact is more marked for the \acrshort{leo} case where the low orbit altitude and the presence of atmospheric drag further slow down the numerical integration.\par

\section{Conclusions}\label{sec:conclusions}

A novel \acrlong{mf} method for \acrlong{up} is presented. The approach leverages \acrlong{da} techniques and \acrlongpl{gmm} to efficiently map the initial uncertainties through arbitrary nonlinear transformations. A \acrshort{da}-based measure of nonlinearity is computed to detect departure from linearity of each kernel and a splitting algorithm is employed to adapt on demand the number of \acrshort{gmm} components needed to satisfy the desired accuracy. Second-order Taylor expansions of the function output with respect to the initial uncertainty are then used to efficiently map the initial statistics using \acrlong{ut} sigma points. The combined \acrshort{loads}-\acrshort{gmm} algorithm is run on a \acrshort{lf} model to minimize the associated computational load. A \acrlong{hf} propagation of the kernels centers is inserted before the \acrshort{ut} sampling to improve the accuracy of the low-fidelity maps. The \acrshort{mf} algorithm is then applied to the problem of orbit \acrshort{up} and three test cases are provided to cover the high-, medium- and low- Earth orbit regimes. The estimated \acrshortpl{pdf} are compared with reference \acrlong{mc} simulations and a good agreement between \acrshort{mc} samples and \acrshort{pdf} contour lines is obtained. The advantages of an adaptive split scheme for the determination of the optimal number of kernels are then highlighted comparing solutions obtained for increasing number of \acrshort{gmm} components and the importance of the \acrshort{hf} correction is assessed comparing the agreement to \acrshort{mc} samples of the \acrshort{lf} and \acrshort{mf} solutions. The proposed \acrshort{mf} approach is estimated to bring a speedup of 15-20 times with respect to its \acrshort{hf} counterpart for the same overall accuracy.

\section*{Acknowledgments}

This work is co-funded by the CNES through A. Fossà’s PhD program, and made use of the CNES orbital propagation tools, including the PACE library.

\bibliography{references}

\end{document}